\input amstex\documentstyle{amsppt}  
\pagewidth{12.5cm}\pageheight{19cm}\magnification\magstep1
\topmatter
\title{Hecke modules based on involutions in extended Weyl groups}\endtitle
\author G. Lusztig\endauthor
\address{Department of Mathematics, M.I.T., Cambridge, MA 02139}\endaddress
\abstract{ Let $X$ be the group of weights of a maximal torus of a simply connected
semisimple group over $\bold C$ and let $W$ be the Weyl group. The semidirect product
$W((\bold Q\otimes X)/X)$ is called an extended Weyl group. There is a natural $\bold C(v)$-algebra
$\bold H$ called the extended Hecke algebra with basis 
indexed by the extended Weyl group which contains the usual Hecke algebra
as a subalgebra. We construct an $\bold H$-module with basis indexed by the involutions in the
extended Weyl group. This generalizes a construction of the author and Vogan.}\endabstract
\thanks{Supported by NSF grant DMS-1566618.}\endthanks
\endtopmatter   
\document

\define\dz{\dot z}

\define\dw{\dot w}

\define\ds{\dot s}

\define\unu{\un{\nu}}

\define\mpb{\medpagebreak}

\define\si{\sim}

\define\sqc{\sqcup}

\define\bX{\bar X}

\define\lb{\linebreak}

\define\op{\oplus}
   
\define\part{\partial}
\define\emp{\emptyset}

\define\ra{\rangle}
\define\rf{\rfloor}
\define\n{\notin}

\define\m{\mapsto}
\define\do{\dots}
\define\la{\langle}
\define\lf{\lfloor}
\define\bsl{\backslash}

\define\sm{\smallmatrix}
\define\esm{\endsmallmatrix}
\define\sub{\subset}    

\define\T{\times}
\define\ti{\tilde}
\define\nl{\newline}
\redefine\i{^{-1}}
\define\fra{\frac}
\define\un{\underline}

\define\ot{\otimes}

\define\Hom{\text{\rm Hom}}

\define\a{\alpha}
\redefine\b{\beta}

\define\g{\gamma}
\redefine\d{\delta}
\define\e{\epsilon}

\define\io{\iota}

\define\p{\pi}
\define\ph{\phi}
\define\ps{\psi}

\define\s{\sigma}
\redefine\t{\tau}
\define\th{\theta}
\define\k{\kappa}
\redefine\l{\lambda}
\define\z{\zeta}
\define\x{\xi}

\redefine\G{\Gamma}
\redefine\D{\Delta}

\define\Ph{\Phi}

\redefine\aa{\bold a}

\define\kk{\bold k}

\define\qq{\bold q}

\define\CC{\bold C}

\define\HH{\bold H}

\define\MM{\bold M}
\define\NN{\bold N}

\define\QQ{\bold Q}
\define\RR{\bold R}

\define\ZZ{\bold Z}

\define\cf{\Cal F}

\define\cl{\Cal L}

\define\cn{\Cal N}
\define\co{\Cal O}

\define\cs{\Cal S}
\define\ct{\Cal T}

\define\cz{\Cal Z}
\define\cx{\Cal X}
\define\cy{\Cal Y}

\define\fe{\frak e}

\define\fs{\frak s}

\define\fH{\frak H}

\define\fP{\frak P}
\define\fQ{\frak Q}

\define\ta{\ti a}

\define\tc{\ti c}

\define\tu{\ti u}

\define\tz{\ti z}

\define\tW{\ti W}
\define\tX{\ti X}

\define\tfQ{\ti{\fQ}}

\define\sha{\sharp}

\define\bul{\bullet}

\define\che{\check}
\define\cha{\che{\a}}

\define\cir{\circ}

\define\CDGVII{L1}
\define\BAR{L2}
\define\MND{L3}
\define\LIF{L4}
\define\LV{LV}
\define\MS{MS}
\define\TIT{T}
\define\YOK{Y}

\head Introduction and statement of results\endhead
\subhead 0.1\endsubhead
Let $\kk$ be an algebraically closed field. Let $G$ be a connected reductive
group over $\kk$. 
Let $T$ be a maximal torus of $G$ and let $U$ be the unipotent radical of a
Borel subgroup of $G$ containing $T$. Let $N$ be the normalizer of $T$
and let $W=N/T$ be the Weyl group; let $w\m|w|$ be the length function on $W$, let 
$S=\{w\in W;|w|=1\}$ and let $\k:N@>>>W$ be the obvious map. The obvious action of $W$ on $T$ 
is denoted by $w:t\m w(t)$.
Let $Y=\Hom(\kk^*,T),X=\Hom(T,\kk^*)$ and let $\la,\ra:Y\T X@>>>\ZZ$ be the obvious pairing. 
We regard $Y,X$ as groups with operation written as addition. Let $K$ be a field of characteristic zero
and let $X_K=K\ot X=\Hom(Y,K)$. Let $\bX=X_K/X=(K/\ZZ)\ot X$. 
The obvious pairing $\la,\ra:Y\T X_K@>>>K$ restricts to a pairing $Y\T X@>>>\ZZ$ hence it induces
a pairing $\lf,\rf:Y\T\bX@>>>K/\ZZ$.
We define an action of $W$ on $Y$ by $w:y\m y'$ where $y'(z)=w(y(z))$ for $z\in\kk^*$. We define an action of $W$ on 
$X_K$ by the equality $\la w(y),w(x)\ra=\la y,x\ra$ for all $y\in Y,x\in X_K,w\in W$.
This action preserves $X$ hence it induces a $W$-action on $\bX$. Let $\che R\sub Y$ be the set of coroots,
let $\che R^+\sub\che R$ be the set of positive coroots determined by $U$, let 
$\che R^-=\che R-\che R^+$. For $s\in S$ we denote by $\cha_s\in Y$ 
the simple coroot such that $s(\cha_s)=-\cha_s$. For $\l\in\bX,s\in S$ we write $s\in W_\l$ if
$\lf\cha_s,\l\rf=0$; we write $s\n W_\l$ if $\lf\cha_s,\l\rf\ne0$. Note that if $s\in W_\l$ then $s\l=\l$.
For $s\in S$ let $T_s$ be the image of $\cha_s:\kk^*@>>>T$.

\subhead 0.2\endsubhead
Let $W_2=\{w\in W;w^2=1\}$. For any integer $m\ge1$ we set 
$$\bX_m=\{\l\in\bX;m^2\l=\l\},$$
$$\tX_m=\{(w,\l)\in W_2\ot\bX;w(\l)=-m\l\}.$$

We write $W\bX$ instead of $W\T\bX$ with the group structure
$(w,\l)(w',\l')=(ww',w'{}\i(\l)+\l')$. We call $W\bX$ the {\it extended Weyl group}. Then
$$\tX_1=\{(w,\l)\in W_2\T\bX;w(\l)=-\l\}=\{(w,\l)\in W\bX;(w,\l)^2=(1,0)\}$$
is exactly the set of involutions in the extended Weyl group $W\bX$.

More generally, if $m\ge1$, then $\{(w,\l)\in W\T\bX;\l\in\bX_m\}$ is a subgroup
of $W\bX$ denoted by $W\bX_m$ and $(w,\l)\m(w,\l)^*:=(w,m\l)$ is an involutive
automorphism of $W\bX_m$. Moreover, $\tX_m$ is the set of $*$-twisted involutions
of $W\bX_m$, that is the set of all $(w,\l)\in W\bX_m$ such that $(w,\l)(w,\l)^*=(1,0)$.

If $m\ge1$ and $(w,\l)\in\tX_m$ then $\l\in\bX_m$.
Note that if $(w,\l)\in\tX_m$ and $s\in S$ then $(sws,s\l)\in\tX_m$;
if in addition $sw=ws$, then $(w,s\l)\in\tX_m$. If we have both $sw=ws$ and
$s\l=\l$ then $(sw,\l)\in\tX_m$. 

\mpb

Let $p$ be a prime number and let $q>1$ be a power of $p$. We set $Q=q^2$. We assume that 
the characteristic of $\kk$ is either $0$ or $p$. Then $\bX_q,\tX_q$ are defined.

We fix a square root $\sqrt{-1}$ of $-1$ in $\CC$.
For $\l\in\bX_q$, $s\in S$, we define $[\l,s]\in\{1,-1\}$ as follows.
We have $\la\cha_s,\l\ra=e/(Q-1)$ with $e\in\ZZ$. When $p\ne2$ we set $[\l,s]=1$ if 
$e\in2\ZZ$ and $[\l,s]=\sqrt{-1}$ if $e\in\ZZ-2\ZZ$; when $p=2$ we set $[\l,s]=1$.

\subhead 0.3\endsubhead
For $w\in W_2,s\in S$ such that $sw=ws$ we define, following \cite{\LIF, 1.18}, a number $(w:s)\in\{-1,0,1\}$ as follows.
Assume first that $G$ is almost simple, simply laced. 
In \cite{\LIF, 1.5, 1.7}, a root system with set of coroots $\che R_w\sub\che R$
and a set of simple coroots $\che\Pi_w$ for $\che R_w$ was associated to $w$; we have $\cha_s\in\che\Pi_w$. This root system
is simply laced and has no component of type $A_l,l>1$. If the component containing $\cha_s$ is not of type $A_1$, there is a unique sequence
$\cha_1,\cha_2,\do,\cha_e$ in $\che\Pi_w$ such that
$\cha_i,\cha_{i+1}$ are joined in the Dynkin diagram of $\che R_w$
for $i=1,2,\do,e-1$,  $\cha_1=\cha_s$  and $\cha_e$ corresponds to a branch point of the Dynkin diagram of $\che R_w$;
if the component containing $\cha_s$ is of type $A_1$ we define $\cha_1,\cha_2,\do,\cha_e$ as the sequence with one term $\cha_s$ (so that
$e=1$). We define $(w:s)=(-1)^e$ if $|sw|<|w|$ and $(w:s)=(-1)^{e+1}$ if $|sw|>|w|$.
Next we assume that $G$ is almost simple, simply connected, not simply laced.
Then $G$ can be regarded as a fixed point set of an automorphism  of a simply connected almost simple, simply laced group $G'$
(as in \cite{\LIF, 1.14}) with Weyl group $W'$, a Coxeter group with a length preserving automorphism $W'@>>>W'$ with
fixed point set $W$. When $s$ is regarded as an element of $W'$, it is a product of $k$ commuting simple reflections
$s'_1,s'_2,\do,s'_k$ of $W'$; here $k\in\{1,2,3\}$. If $k\ne2$, we define $(w:s)$ for $W$ to be $(w:s_i)$ for $G'$ where
$i$ is any element of $\{1,\do,k\}$. If $k=2$ we have either $ws_1=s_1w$, $ws_2=s_2w$
(and $(w:s)$ for $G$ is defined to be $(w:s_1)=(w:s_2)$ for $G'$) or $ws_1=s_2w$, $ws_2=s_1w$ (and $(w:s)$ for $G$ is defined to be $0$.)
We now drop the assumption that $G$ is almost simple. Let $G''$ be the simply connected cover of an
almost simple factor of the adjoint group of $G$ with Weyl group $W''\sub W$
such that $s\in W''$ and let $w''$ be the $W'$-component of $w$. Then $(w:s)$ for $G$ is defined to be $(w'':s)$ for $G''$ (which is
is defined as above).

For $p,q$ as in 0.2, $(w,\l)\in\tX_q$, $s\in S$ such that $sw=ws$, we set
$$\d_{w,\l;s}=\exp(2\pi\sqrt{-1}((q-e)/2)(1-(w:s))\la\cha_s,\l\ra)$$
if $p\ne2$,  $e=|w|-|sw|=\pm1$ and $\d_{w,\l;s}=1$ if $p=2$. (Note that
$\exp(2\pi\sqrt{-1}x)$ is well defined for $x\in\QQ/\ZZ$.)
If $G$ is simply laced, then $\d_{w,\l;s}=1$ (since $(w:s)=\pm1$).
In general we have $\d_{w,\l;s}=\pm1$. Indeed, we can assume that $p\ne2$.
It is enough to show that $(q-e)\la\cha_s,\l\ra=0$. From our assumption we have
$$\lf\cha_s,\l\rf=\lf w\cha_s,w\l\rf=\lf-e\cha_s,-q\l\rf=qe\lf\cha_s,\l\rf=qe\i\lf\cha_s,\l\rf$$
hence $(q-e)\lf\cha_s,\l\rf=0$ and our claim follows.

\mpb

The following assumption will be made in parts of the paper (it will simplify some proofs).

(a) For $s\in S$, $\cha_s;\kk^*@>>>T_s$ is an isomorphism.
\nl
This is certainly satisfied if $G$ is simply connected. 

Here is one of the main results of this paper:
\proclaim{Theorem 0.4}Let $q,p$ be as in 0.2. Assume that 0.3(a) holds.
Let $M_q$ be the $\CC$-vector space with basis $\{a_{w,\l};(w,\l)\in\tX_q\}$.
If $p\ne2$ let $z\in\ZZ$ be such that $2z\n(q^2-1)\ZZ$; if $p=2$ let $z\in\ZZ$ be arbitrary. There is a unique action of 
the braid group of $W$ on $M_q$ in which the generators $\{\ct_s;s\in S\}$ of the braid group applied to the basis elements
of $M_q$ are as follows. (We set $\D=1$ if $s\in W_\l$ and $\D=0$ if $s\n W_\l$.)

(a) $\ct_sa_{w,\l}=a_{sws,\l}$ if $sw\ne ws,|sw|>|w|,\D=1$;

(b) $\ct_sa_{w,\l}=a_{sws,\l}+(q-q\i)a_{w,\l}$ if $sw\ne ws,|sw|<|w|,\D=1$;

(c) $\ct_sa_{w,\l}=a_{w,\l}+(q+1)a_{sw,\l}$ if $sw=ws,|sw|>|w|,\D=1$;

(d) $\ct_sa_{w,\l}=(1-q\i)a_{sw,\l}+(q-q\i-1)a_{w,\l}$ if $sw=ws,|sw|<|w|,\D=1$;

(e)  $\ct_sa_{w,\l}=[\l,s]a_{sws,s\l}$ if $sw\ne ws,|sw|>|w|,\D=0$;

(f) $\ct_sa_{w,\l}=[\l,s]\i a_{sws,s\l}$ if $sw\ne ws,|sw|<|w|,\D=0$;

(g) $\ct_sa_{w,\l}=\d_{w,s\l;s}a_{w,s\l}$ if $sw=ws, |sw|>|w|,\D=0$;

(h) $\ct_sa_{w,\l}=-\d_{w,s\l;s}\exp(2\pi\sqrt{-1}(w:s)z\la\cha_s,\l\ra)a_{w,s\l}$ 
if $sw=ws,|sw|<|w|,\D=1$.
\endproclaim
Note that the subspace of $M_q$ spanned by $\{a_{w,0};w\in W_2\}$ is stable under the
braid group action; the resulting braid group action on that subspace involves only
the cases where $\D=1$ and in fact is the representation of the Hecke algebra of $W$ with parameter $q$ introduced
in \cite{\LV}.
Thus the theorem is a generalization of a part of \cite{\LV}. 
In the general case we can define operators $1_\l:M_q@>>>M_q$ (for $\l\in\bX_q$) by
$1_\l a_{w,\l'}=\d_{\l,\l'}a_{w,\l'}$ for all $(w,\l')\in\tX_q$.
The operators $\ct_s$ and $1_\l$ on $M_q$ satisfy the relations
of an "extended Hecke algebra", isomorphic to the endomorphism algebra of the representation of 
$G(F_q)$ induced by the trivial representation of $U(F_q)$ (assuming that $\kk$ is an algebraic closure of a finite
field $F_q$ and $G$ is split over $F_q$). This endomorphism algebra was studied by
Yokonuma \cite{\YOK} and a description of it in terms of generators like $\ct_s,1_\l$
was given in \cite{\CDGVII}. 
The proof of the theorem is given in \S4, in terms of $G(F_q),U(F_q)$ as above.
Namely we show that $M_q$ can be interpreted as the vector space spanned by the double cosets
$\G_1\bsl\G/\G_2$ regarded naturally as a module over the algebra spanned as a vector space by the
 double cosets $\G_1\bsl\G/\G_1$
for suitable finite groups $\G_1\sub\G\supset\G_2$. (In our case we have
$\G=G(F_{q^2})$, $\G_1=U(F_{q^2})$, $\G_2=G(F_q)$.) A key role in our proof is
played by a certain non-standard lifting (introduced in \cite{\LIF}) to $N$ for the involutions in $W$.
(The usual lifting, due to Tits \cite{\TIT}, is not suitable for the purposes of this paper.)

\subhead 0.5\endsubhead
We now assume that $\kk=\CC$. Let $v$ be an indeterminate and let $\MM$
be the $\CC(v)$-vector space with basis $\{\aa_{w,\l};(w,\l)\in\tX_1\}$.
For any

(a) $(w,\l)\in\tX_1$ and $s\in S$ such that $|sw|>|w|$
\nl
we set
$$\d'_{w,\l;s}=\exp(2\pi\sqrt{-1}(1-(w:s))\lf\cha_s,\l\rf).$$
 We note that for $w,\l,s$ as in (a) we have
$$\lf\cha_s,\l\rf=\lf w\cha_s,w\l\rf=\lf\cha_s,-\l\rf=-\lf\cha_s,\l\rf$$
hence

(b) $2\lf\cha_s,\l\rf=0$
\nl
so that $\d'_{w,\l;s}$ is well defined and is in $\{1,-1\}$.
The following result is a generic version of Theorem 0.4 in which
$q$ is replaced by $v^2$  and $M_q$ is replaced by $\MM$.

\proclaim{Theorem 0.6} We assume that $\kk=\CC$ and that 0.3(a) holds.
There is a unique action of the braid group of $W$ on
$\MM$ in which the generators $\{\ct_s;s\in S\}$ of the braid group applied to
the basis elements of $\MM$ are as follows. (We write $\D=1$ if $s\in W_\l$
and $\D=0$ if $s\n W_\l$.)

(a) $\ct_s\aa_{w,\l}=\aa_{sws,s\l}$ if $sw\ne ws,|sw|>|w|$;

(b) $\ct_s\aa_{w,\l}=\aa_{sws,s\l}+\D(v^2-v^{-2})\aa_{w,\l}$ if
$sw\ne ws,|sw|<|w|$;

(c) $\ct_s\aa_{w,\l}=\d'_{w,s\l;s}\aa_{w,s\l}+\D(v+v\i)\aa_{sw,\l}$ if $sw=ws,|sw|>|w|$;

(d) $\ct_s\aa_{w,\l}=\D(v-v\i)\aa_{sw,\l}+\D(v^2-v^{-2})\aa_{w,\l}-\aa_{w,s\l}$ if $sw=ws,|sw|<|w|$.
\endproclaim
This can be deduced from Theorem 0.4 (see \S4).

We can interpret the theorem as providing an $\HH$-module structure on $\MM$ where
$\HH$ is the extended Hecke algebra (see 4.5).
The subspace of $\MM$ spanned by $\{\aa_{w,0};w\in W_2\}$ is stable
under the operators $\ct_s$ and this defines a representation of
the generic Hecke algebra of $W$ which was defined in \cite{\LV}.

\subhead 0.7\endsubhead
The  action in 0.6 can be specialized to $v=1$. It becomes the braid group
action on the $\CC$-vector space with basis $\{\aa_{w,\l};(w,\l)\in\tX_1\}$ in
which the generators $\ct_s$ of the braid group act as follows.
(Notation and assumptions of 0.6.)

(a) $\ct_s\aa_{w,\l}=\aa_{sws,s\l}$ if $sw\ne ws$;

(b) $\ct_s\aa_{w,\l}=\d'_{w,s\l;s}\aa_{w,s\l}+2\D\aa_{sw,\l}$ if $sw=ws,|sw|>|w|$;

(c) $\ct_s\aa_{w,\l}=-\aa_{w,s\l}$ if $sw=ws,|sw|<|w|$.
\nl
This is actually a $W$-action.

\subhead 0.8\endsubhead
Let $m$ be an integer $\ge1$ and let $\MM_m$ be the $\CC(v)$-vector space with basis $\{\aa_{w,\l};(w,\l)\in\tX_m\}$.
In the following result (a variant of 0.4 and 0.6) the assumption 0.3(a) is not used.

\proclaim{Theorem 0.9} There is a unique action of the braid group of $W$ on
$\MM_m$ in which the generators $\{\ct_s;s\in S\}$ of the braid group applied to
the basis elements of $\MM_m$ are as follows. (We write $\D=1$ if $s\in W_\l$
and $\D=0$ if $s\n W_\l$.)

(a) $\ct_s\aa_{w,\l}=\aa_{sws,s\l}$ if $sw\ne ws,|sw|>|w|$;

(b) $\ct_s\aa_{w,\l}=\aa_{sws,s\l}+\D(v^2-v^{-2})\aa_{w,\l}$ if
$sw\ne ws,|sw|<|w|$;

(c) $\ct_s\aa_{w,\l}=\aa_{w,s\l}+\D(v+v\i)\aa_{sw,\l}$ if $sw=ws,|sw|>|w|$;

(d) $\ct_s\aa_{w,\l}=\D(v-v\i)\aa_{sw,\l}+\D(v^2-v^{-2}-1)
\aa_{w,\l}+(1-\D)\aa_{w,s\l}$ if $sw=ws,|sw|<|w|$.
\endproclaim
The proof is given in \S6. It relies on results in \cite{\LV} and \cite{\MND}.

\subhead 0.10\endsubhead
The  action in 0.9 can be specialized to $v=1$. It becomes the braid group
action on the $\CC$-vector space with basis $\{\aa_{w,\l};(w,\l)\in\tX_m\}$ in
which the generators $\ct_s$ of the braid group act as follows.
(Notation and assumptions of 0.9.)

(a) $\ct_s\aa_{w,\l}=\aa_{sws,s\l}$ if $sw\ne ws$;

(b) $\ct_s\aa_{w,\l}=\aa_{w,s\l}+2\D\aa_{sw,\l}$ if $sw=ws,|sw|>|w|$;

(c) $\ct_s\aa_{w,\l}=\aa_{w,s\l} -2\D\aa_{w,\l}$ if $sw=ws,|sw|<|w|$.
\nl
This is actually a $W$-action.

\subhead 0.11\endsubhead
{\it Notation.} If $X\sub X'$ are sets and $\io:X'@>>>X'$ satisfies $\io(X)\sub X$ we write $X^\io=\{x\in X;\io(x)=x\}$.

\subhead 0.12\endsubhead
I thank David Vogan for discussions.

\head 1. The algebra $\cf$\endhead
\subhead 1.1\endsubhead
Let $p,q,Q$ be as in 0.2. We now assume that $\kk$ is an algebraic closure of the finite
field $F_q$ with $\sha(F_q)=q$. We fix a pinning $(x_s:\kk@>>>G,y_s:\kk@>>>G;s\in S)$ corresponding 
to $T,U$. (We have $x_s(\kk)\sub U$.) 
Let $W@>>>N$, $w\m\dw$ be the Tits cross-section of $\k:N@>>>W$ associated to this pinning, see \cite{\TIT}.
We fix an $F_q$-rational structure on $G$ with Frobenius map $\ph:G@>>>G$ such that
$\ph(t)=t^q$ for all $t\in T$ and $\ph(x_s(z))=x_s(z^q),\ph(y_s(z))=y_s(z^q)$ for all
$z\in\kk$. We have $\ph(\dw)=\dw$ for any $w\in W$ and $\ph(U)=U$.
Let $F_Q$ be the subfield of $\kk$ with $\sha(F_Q)=Q$. We set $\Ph=\ph^2$. We set $\e=-1\in\kk^*$.

For $s\in S$, $z\in\kk^*$ we set $z_s=\cha_s(z)\in T_s$.
In particular, $\e_s\in T_s$ is defined and we have $\ds^2=\e_s$.

\subhead 1.2\endsubhead
Let $\cx=G/U$. Now $G$ acts on $\cx$ by $g:xU\m gxU$ and on $\cx^2$ by 
$g:(xU,yU)\m(gxU,gyU)$. We have $\cx^2=\sqc_{n\in N}O_n$ where $O_n=\{(xU,yU)\in\cx^2;x\i y\in UnU\}$.
Now $\ph,\Ph$ induce endomorphisms of $\cx$ and $\cx^2$ denoted again by $\ph,\Ph$. 
For $n\in N$, we have $\ph(O_n)=O_{\ph(n)}$ hence $\Ph(O_n)=O_{\Ph(n)}$. Thus we have
$(\cx^2)^\Ph=\sqc_{n\in N^\Ph}O_n^\Ph$ and
$O_n^\Ph$ ($n\in N^\Ph$) are exactly the orbits of $G^\Ph$ on $(\cx^2)^\Ph$.

\subhead 1.3\endsubhead
Let
$$\cf=\{f:(\cx^2)^\Ph@>>>\CC;f\text{ constant on the orbits of }G^\Ph\}.$$
This is a $\CC$-vector space with basis $\{k_n;n\in N^\Ph\}$ where $k_n$ is $1$ on
$O_n^\Ph$ and is $0$ on $(\cx^2)^\Ph-O_n^\Ph$.
Now $\cf$ is an associative algebra with $1$ under convolution:
$$(f_1f_2)(xU,zU)=\sum_{yU\in\cx^\Ph}f_1(xU,yU)f_2(yU,zU);$$
 here $f_1\in\cf,f_2\in\cf$, $(xU,zU)\in(\cx^2)^\Ph$. 

We now state two well known results. 

\proclaim{Lemma 1.4} Assume that $n,n'\in N$, $\k(n)=w,\k(n')=w'$ satisfy $|ww'|=|w|+|w'|$. 

(a) If $(xU,yU)\in O_n,(yU,zU)\in O_{n'}$ then $(xU,zU)\in O_{nn'}$.

(b) If $(xU,zU)\in O_{nn'}$ then there is a unique $yU\in X$ such that
$(xU,yU)\in O_n,(yU,zU)\in O_{n'}$.
\endproclaim

\proclaim{Lemma 1.5} Assume that $s\in S$. Assume that 0.3(a) holds.

(a) If $(xU,x'U)\in O_{\ds},(x'U,zU)\in O_{\ds\i}$ then $(xU,zU)\in O_1$ or
$(xU,zU)\in\sqc_{y\in T_s}O_{\ds y}$.

(b) If $(xU,zU)\in O_1$ then 
$\{x'U\in X;(xU,x'U)\in O_{\ds},(x'U,zU)\in O_{\ds\i}\}$ is an affine line.

(c) If $(xU,zU)\in O_{\ds y}$ with $y\in T_s$, then 
$\{x'U\in X;(xU,x'U)\in O_{\ds},(x'U,zU)\in O_{\ds\i}\}$ is a point.
\endproclaim

The following result can be deduced from Lemmas 1.4, 1.5.
\proclaim{Lemma 1.6} Assume that $s\in S$, $n\in N$, $\k(n)=w$ satisfy $|ws|<|w|$.
Assume that 0.3(a) holds.

(a) If $(xU,x'U)\in O_n,(x'U,x''U)\in O_{\ds\i}$ then $(xU,x''U)\in O_{n\ds\i}$ or \lb
$(xU,x''U)\in\sqc_{\t\in T_s}O_{n\t}$.

(b) If $(xU,x''U)\in O_{n\ds\i}$ then 
$\{x'U\in X;(xU,x'U)\in O_n,(x'U,x''U)\in O_{\ds\i}\}$ is an affine line.

(c) If $(xU,x''U)\in O_{n\t}$ with $y\in T_s$, then 
$$\{x'U\in X;(xU,x'U)\in O_n,(x'U,x''U)\in O_{\ds\i}\}$$
 is a point.
\endproclaim

\subhead 1.7\endsubhead
Assume that 0.3(a) holds.
From 1.4 we deduce that for $n,n'\in N^\Ph$ such that $|\k(nn')|=|\k(n)|+|\k(n')|$ 
we have
$$k_nk_{n'}=k_{nn'},\tag a$$
in $\cf$. In particular $k_1$ is the unit element of $\cf$. From 1.5 we deduce that for $s\in S$ we have
$$k_{\ds} k_{\ds}=Qk_{\e_s}+\sum_{y\in T_s^\Ph}k_{\ds} k_y.\tag b$$
(See \cite{\YOK}.) It follows that for $s\in S,w\in W,n\in N^\Ph$ such that $|sw|<|w|$, $\k(n)=w$ we have
$$k_{\ds}k_n=Qk_{\ds n}+\sum_{y\in T_s^\Ph}k_{yn}\tag c$$
and for $s\in S,w\in W,n\in N^\Ph$ such that $|ws|<|w|$, $\k(n)=w$ we have
$$k_nk_{\ds\i}=Qk_{n\ds\i}+\sum_{y\in T_s^\Ph}k_{ny}.\tag d$$
From (a),(c),(d) we deduce that for $s\in S,w\in W,n\in N^\Ph$ such that $sw=ws,|sw|<|w|$, $\k(n)=w$ we have
$$k_{\ds} k_nk_{\ds\i}=Qk_{\ds n\ds\i}+Q\sum_{y\in T_s^\Ph}k_{\ds ny}+\sum_{y\in T_s^\Ph,y'\in T_s^\Ph}k_{yny'}.
\tag e$$

\subhead 1.8\endsubhead
We set $\fs=\Hom(T^\Ph,\CC^*)$. Now $W$ acts on $\fs$ by
$w:\nu\m w\nu$ where $(w\nu)(t)=\nu(w\i(t))$ for $t\in T^{\Ph^h}$. For $\nu\in\fs$ we set 
$$1_\nu=|T^\Ph|\i\sum_{\t\in T^\Ph}\nu(\t)k_\t\in\cf.\tag a$$ 
We have
$$\sum_{\nu\in\fs}1_\nu=k_1=1.\tag b$$
Indeed,
$$\sum_{\nu\in\fs}1_\nu=|T^\Ph|\i\sum_{\t\in T^\Ph}\sum_{\nu\in\fs}\nu(\t)k_\t
=\sum_{\t\in T^\Ph}\d_{\t,1}k_\t=k_1.$$
For $\nu,\nu'$ in $\fs$ we have
$$1_\nu1_{\nu'}=\d_{\nu,\nu'}1_\nu.\tag c$$
Indeed,
$$\align&1_\nu1_{\nu'}=|T^\Ph|^{-2}\sum_{\t\in T^\Ph,\t'\in T^\Ph}\nu(\t)\nu'(\t')k_{\t\t'}\\&=
|T^\Ph|^{-2}\sum_{\t\in T^\Ph,\t''\in T^\Ph}\nu(\t)\nu'(\t''\t\i))
k_{\t''}\\&=\d_{\nu,\nu'}|T^\Ph|\i\sum_{\t''\in T^\Ph}\nu'(\t'')k_{\t''}=\d_{\nu,\nu'}1_\nu.\endalign$$
For $\nu\in\fs,n\in N^\Ph,w=\k(n)\in W$ we have
$$k_n1_\nu=1_{w\nu}1_\nu.\tag d$$
Indeed,
$$\align&k_n1_\nu=|T^\Ph|\i\sum_{\t\in T^\Ph}\nu(\t)k_{n\t}=|T^\Ph|\i\sum_{\t\in T^\Ph}\nu(\t)k_{w(\t)n}\\&=
|T^\Ph|\i\sum_{\t'\in T^\Ph}\nu(w\i(\t'))k_{\t'n}=1_{w\nu}k_n.\endalign$$
For $t\in T^\Ph,\nu\in\fs$ we have 
$$k_t1_\nu=\nu(t\i)1_\nu.\tag e$$
Indeed,
$$\align&k_t1_\nu=|T^\Ph|\i\sum_{\t\in T^\Ph}\nu(\t)k_{t\t}=|T^\Ph|\i\sum_{\t'\in T^\Ph}\nu(t\i\t')k_{\t'}\\&=
\nu(t\i)|T^\Ph|\i\sum_{\t'\in T^\Ph}\nu(\t')k_{\t'}=\nu(t\i)1_\nu.\endalign$$
For $\nu\in\fs,s\in S$ we write $s\in W_\nu$ if $\nu(\cha_s(z))=1$ for all $z\in F_Q^*$ or equivalently
if $\nu|_{T_s^\Ph}=1$; we write $s\n W_\nu$ if $\nu|_{T_s^\Ph}$ is not identically $1$.

For $\nu\in\fs,\cha\in\che R$ we define $[\nu,\cha]$ as follows.
If $\nu(\cha(\e))=1$ we set $[\nu,\cha]=1$; if $\nu(\cha(\e))=-1$ we set $[\nu,\cha]=\sqrt{-1}$.
(Since $\cha(\e)^2=1$ we must have $\nu(\cha(\e))\in\{1,-1\}$.) If $p=2$ we have $\cha(\e)=1$ hence
$[\nu,\cha]=1$. We have $[\nu,\cha]^2=\nu(\cha(\e))$.

For $s\in S$ we set
$$\ct_s=q\i k_{\ds}\sum_{\nu\in\fs}[\nu,\cha_s]1_\nu\in\cf.\tag f$$  
We show:
$$\ct_s\ct_s=1+(q-q\i)\sum_{\nu\in\fs;s\in W_\nu}\ct_s1_\nu.\tag g$$
Indeed, we have
$$\align&\ct_s\ct_s=Q\i\sum_{\nu\in\fs,\nu'\in\fs}[\nu,\cha_s][\nu',\cha_s]k_{\ds}1_\nu k_{\ds}1_{\nu'}\\&=
Q\i\sum_{\nu\in\fs,\nu'\in\fs}[\nu,\cha_s][\nu',\cha_s]k_{\ds}k_{\ds}1_{s\nu}1_{\nu'}\\&=
Q\i\sum_{\nu'\in\fs}[s\nu',\cha_s][\nu',\cha_s]k_{\ds}k_{\ds}1_{\nu'}\\&=
Q\i\sum_{\nu\in\fs}\nu(\e_s)k_{\ds}k_{\ds}1_\nu\\&
=\sum_{\nu\in\fs}\nu(\e_s)k_{\e_s}1_\nu+Q\i\sum_{\nu\in\fs,y\in T_s^\Ph}\nu(\e_s)k_{\ds} k_y1_\nu\\&
=\sum_{\nu\in\fs}1_\nu+Q\i\sum_{\nu\in\fs,y\in T_s^\Ph}\nu(\e_s)\nu(y\i)k_{\ds}1_\nu\\&
=1+Q\i(Q-1)\sum_{\nu\in\fs,\nu|_{T_s^\Ph}=1}k_{\ds}1_\nu.\endalign$$
It remains to use that if $\nu|_{T_s^\Ph}=1$ then $\nu(\e_s)=1$ hence $[\nu,\cha_s]=1$.

Now (g) implies that $\ct_s\i\in\cf$ is well defined and we have
$$\ct_s\i=\ct_s-(q-q\i)\sum_{\nu\in\fs;s\in W_\nu}1_\nu.\tag h$$
From (h) we see that for any $\nu\in\fs$:
$$\ct_s\i1_\nu=\ct_s1_\nu-\D(q-q\i)1_\nu\tag i$$
where $\D=1$ if $s\in W_\nu$, $\D=0$ if $s\n W_\nu$.

For any $\nu\in\fs$ we show:
$$1_\nu\ct_s=\ct_s1_{s\nu}.\tag j$$   
Indeed, we have
$$1_\nu\ct_s=q\i1_\nu k_{\ds}\sum_{\nu'\in\fs}[\nu',\cha_s]1_{\nu'}=q\i\sum_{\nu'\in\fs}k_{\ds} [\nu',\cha_s]1_{s\nu}1_{\nu'}=
q\i k_{\ds} [\nu,\cha_s]1_{s\nu},$$
$$\ct_s1_{s\nu}=q\i k_{\ds}\sum_{\nu'\in\fs}[\nu',\cha_s]1_{\nu'}1_{\s\nu}=q\i k_{\ds}[\nu,\cha_s]1_{\s\nu}.$$

\subhead 1.9\endsubhead
For any $w\in W$ we set
$$\ct_w=q^{-|w|}k_{\dw}\sum_{\nu\in\fs}\prod_{\cha\in\che R^+;w\i(\cha)\in\che R^-}[\nu,w\i\cha]1_\nu\in\cf.\tag a$$  
When $w=s\in S$, this definition agrees with the earlier definition of $\ct_s$. For $s\in S$, $w\in W$ such that $|ws|>|w|$ we show
$$\ct_{ws}=\ct_w\ct_s.\tag a$$
Since $|ws|>|w|$, we have $w(\cha_s)\in R^+$ and
$\{\cha\in\che R^+;(ws)\i(\cha)\in\che R^-\}=\{\cha\in R^+;w\i(\cha)\in\che R^-\}\sqc\{w(\cha_s)\}$.
Hence we have
$$\align&\ct_{ws}=q^{-|ws|}
k_{\dw\ds}\sum_{\nu\in\fs}\prod_{\cha\in\che R^+,(ws)\i)(\cha)\in\che R^-}[\nu,(ws)\i\cha]1_{\nu}\\&=
q^{-|ws|}k_{\dw\ds}\sum_{\nu\in\fs}[\nu,(ws)\i(w(\cha_s))]\prod_{\cha\in\che R^+;w\i(\cha)\in\che R^-}
[\nu,(ws)\i(\cha)]1_\nu\\&=q^{-|ws|}k_{\dw\ds}\sum_{\nu\in\fs}[\nu,\cha_s]\prod_{\cha\in\che R^+;w\i(\cha)\in\che R^-}
[\nu,(ws)\i(\cha)]1_\nu.\endalign$$
We have
$$\align&\ct_w\ct_s=q^{-|w|}q\i k_{\dw}\sum_{\nu\in\fs,\nu'\in\fs}
\prod_{\cha\in\che R^+,w\i(\cha)\in\che R^-}[\nu,w\i(\cha)][\nu',\cha_s]1_\nu k_{\ds}1_{\nu'}=\\&
q^{-|ws|}k_{\dw}k_{\ds}\sum_{\nu\in\fs,\nu'\in\fs}[\nu',\cha_s]\prod_{\cha\in\che R^+,w\i(\cha)\in\che R^-}[\nu,w\i(\cha)]
1_{s\nu}1_{\nu'}\\&=
q^{-|ws|}k_{\dw\ds}\sum_{\nu\in\fs}[\nu,\cha_s]\prod_{\cha\in\che R^+,w\i(\cha)\in\che R^-}[\nu,(ws)\i(\cha)]1_\nu.\endalign$$
This proves (a).

From (a) we deduce:
$$\ct_{ww'}=\ct_w\ct_{w'}\text{ if $w,w'$ in $W$ satisfy } |ww'|=|w|+|w'|.\tag b$$
Using 1.8(j) and (a) we see that
$$1_\nu\ct_w=\ct_w1_{w\i\nu}\text{ for }w\in W,\nu\in\fs.\tag c$$   
We note that

(d) {\it $\{\ct_w1_\nu;w\in W,\nu\in\fs\}$ is a $\CC$-basis of $\cf$.}
\nl
This follows from the fact that (up to a nonzero scalar) $\ct_w1_\nu$ is equal to 
$$\sum_{\t\in T^\Ph}\nu(\t)k_{\dw\t}.$$

\head 2. The $\cf$-module $\cf'$\endhead
\subhead 2.1\endsubhead
In this section we assume that 1.3(a) holds.
We preserve the setup of 1.1. We define $\ph':N@>>>N$ by $\ph'(n)=\ph(n)\i$. 
We define $\ps:\cx^2@>>>\cx^2$ by $\ps(xU,yU)=(\ph(y)U,\ph(x)U)$. This is a Frobenius map for an
$F_q$-rational structure on $\cx^2$. The $G$-action on $\cx^2$ in 1.2 is compatible with this $F_q$-rational
structure on $\cx^2$ and with the $F_q$-rational structure on $G$ given by $\ph$. It follows that any
$G$-orbit $O_n$ on $\cx^2$ such that $\ps(O_n)=O_n$ satisfies the condition that $O_n^\ps\ne\emp$ and that
$G^\ph$ acts transitively on $O_n^\ps$. (We use Lang's theorem and the connectedness of the stabilizers of the
$G$-action on $O_n$.) For $n\in N$ we have $\ps(O_n)=O_{\ph'(n)}$; thus $\ps(O_n)=O_n$ precisely when $n\in N^{\ph'}$. 
Thus we have $(\cx^2)^\ps=\sqc_{n\in N^{\ph'}}O_n^\ps$ and $O_n^\ps$ (for various $n\in N^{\ph'}$) are precisely the
$G^\ph$-orbits in $(\cx^2)^\ps$. Let 
$$\cf'=\{h:(\cx^2)^\ps@>>>\CC;h\text{ constant on the orbits of }G^\ph\}.$$
This is a $\CC$-vector space with basis $\{\th_m;m\in N^{\ph'}\}$ where $\th_m$ is $1$ on
$O_m^\ps$ and is $0$ on $(\cx^2)^\ps-O_m^\ps$. Now $\cf'$ is a $\cf$-module under convolution
$$(fh)(xU,\ph(x)U)=\sum_{yU\in\cx^\Ph}f(xU,yU)h(yU,\ph(y)U);$$
here $f\in\cf,h\in\cf',(xU,\ph(x)U)\in(\cx^2)^\ps$. (In this $\cf$-module, multiplication by
the unit element of $\cf$ is the identity map of $\cf'$.)

\subhead 2.2\endsubhead
Now $\ph':N@>>>N$ is an $F_q$-structure on $N$ not necessarily compatible with the group structure of $N$.
But it is compatible with the $T\T T$-action on $N$ 
given by $(t_1,t_2):n\m t_1nt_2\i$ and the $F_q$-rational structure on $T\T T$ with Frobenius map
$(t_1,t_2)\m(\ph(t_2),\ph(t_1))$. Hence any $T\T T$-orbit of the action on $N$ which is stable under
$\ph':N@>>>N$ must have a $\ph'$-fixed point. Such an orbit is of the form $\k\i(w)$ with $w\in W$
satisfying $w\i=w$ that is $w\in W_2$. Using Lang's theorem and the connectedness  of the stabilizers of the
$T\T T$-action on $\k\i(w)$, we see that for $w\in W_2$, $\k\i(w)\cap N^{\ph'}$ is nonempty and is exactly one
orbit for the subgroup $\{(t_1,t_2)\in T\T T;(t_1,t_2)=(\ph(t_2),\ph(t_1))\}$ of $T\T T$. Thus,

(a) {\it $N^{\ph'}=\sqc_{w\in W_2}N(w)$ where for any $w\in W_2$, $N(w):=\k\i(w)\cap N^{\ph'}$ is nonempty and is a 
single orbit for the action of $T^\Ph$ on $N^{\ph'}$ given by $t:n\m tn\ph(t)\i$.}
\nl
For $w\in W_2$ we have $N(w)=\{\dw t;t\in T,w(t^q)t\dw^2=1\}$. Let $T(w)=\{t\in T;w(t^q)t=1\}$. Clearly, 

(b) {\it $N(w)$ is a single orbit under right translation by $T(w)$.}
\nl
From the definitions we see that:

(c) {\it For $w\in W_2$, $z\in W$ we have $\dz N(w)\dz\i=N(zwz\i)$.}
\nl
It is enough to show that $\dz N^{\ph'}\dz\i=N^{\ph'}$.
More generally, if $n\in N^\Ph$, then $nN^{\ph'}\ph(n)\i=N^{\ph'}$. This is easily verified.

\mpb

For $w\in W_2$, we define a homomorphism $e_w:T^\Ph@>>>T(w)$ by $\t\m w(\t)\t^{-q}$. We show:

(d) {\it $e_w$ is surjective.}
\nl
Let $t\in T(w)$. By Lang's theorem we have $t=w(\t)\t^{-q}$ for some $\t\in T$. Since $t\in T(w)$ we have 
automatically $\t\in T^\Ph$ and (d) follows. 

\mpb

For $w\in I$, $s\in S$ such that $sw=ws$ we show:

(e) {\it If $|sw|>|w|$ then $\{c_s;c\in F_Q,c^{q+1}=1\}\sub T(w)$; if
$|sw|<|w|$ then $\{c_s;c\in F_Q,c^{q-1}=1\}\sub T(w)$.}
\nl
Assume first that $|sw|>|w|$ and that $c^{q+1}=1$. We have $w(c_s)=c_s$ hence $w(c_s^q)c_s=c_s^{q+1}=1$. Next we 
assume that $|sw|<|w|$ and that $c^{q-1}=1$. We have $w(c_s)=c_s\i$ hence $w(c_s^q)c_s=c_s^{-q+1}=1$. This proves 
(e).

\subhead 2.3\endsubhead
For $n\in N^\Ph,m\in N^{\ph'}$ we have $k_n\th_m=\sum_{m'\in N_*}\cn_{n,m,m'}\th_{m'}$ where 
$$\cn_{n,m,m'}=\sha\{yU\in X^\Ph;(xU,yU)\in O_n^\Ph,(yU,\ph(y)U)\in O_m^\ps\}.$$
We have also
$$\cn_{n,m,m'}=\sha Z_{xU,\ph(x)U}^\ps$$
where 
$$Z_{xU,\ph(x)U}=\{(yU,y'U)\in O_m;(xU,yU)\in O_n,(y'U,\ph(x)U)\in O_{\ph(n)\i}\}$$
with $(xU,\ph(x)U)$ fixed in $O_{m'}^\ps$ (note that $Z_{xU,\ph(x)U}$ is $\ps$-stable).

\proclaim{Lemma 2.4} Assume that $n=t\in T^\Ph$, $m\in N^{\ph'}$. We have
$k_t\th_m=\th_{tm\ph(t)\i}$.
\endproclaim
If $m'\in N^{\ph'}$ satisfies $\cn_{n,m,m'}\ne0$ then from Lemma 1.4 (applied twice) we see that
$Z_{xU,\ph(x)U}$ is a point and $m'=tm\ph(t)\i$; moreover we have $\cn_{n,m,m'}=1$. The result follows.

\proclaim{Lemma 2.5} Assume that $s\in S$, $w\in I$, $m\in N(w)$, $sw\ne ws$, $|ws|>|w|$. 
Recall that $\ds m\ds\i\in N(sws)$. We have 
$$k_{\ds}\th_m=\th_{\ds m\ds\i}.$$
\endproclaim
In this case we have $|sws|=|w|+2$.
If $m'\in N^{\ph'}$ satisfies $\cn_{n,m,m'}\ne0$ then from Lemma 1.4 (applied twice) we see that
$Z_{xU,\ph(x)U}$ (in 2.3 with $n=\ds$) is a point and $m'=\ds m\ph(\ds)\i$; moreover we have $\cn_{n,m,m'}=1$. The
result follows.

\proclaim{Lemma 2.6} Assume that $s\in S$, $w\in I$, $m\in N(w)$, $sw=ws$, $|ws|>|w|$. Write $m=\dw t$ where $t\in T$ 
satisfies $w(t^q)t\dw^2=1$.

(a) We have $\dw s(t)=\ds m\ds\i\in N(w)$. We have $s(t)\i t\e_s=\ds m\i\ds m\in T_s$,
$(\ds m\i\ds m)^{q+1}=1$.

(b) For $y\in T_s$ we have $\ds\dw ty=\ds my\in N(sw)$ if and only if $y^{q-1}=s(t)\i t\e_s=\ds m\i\ds m$.
There are exactly $q-1$ such $y$; they are all automatically in $T_s^\Ph$.

(c) We have
$$k_{\ds}\th_m=q\th_{\ds m\ds\i}+\sum_{y\in T_s;y^{q-1}=\ds m\i\ds m}\th_{\ds my}.$$
\endproclaim
The equalities in (a) are easily checked; the inclusion $\ds n\ds\i\in N(w)$ follows from 2.2(c).
We have $s(t)\i t\e_s\in T_s$. To prove (a) it remains to show that $(s(t)\i t\e_s)^{q+1}=1$. 
We have $\ds\dw^2=\dw^2\ds$ hence $\dw^2=\ds\dw^2\ds\i=s(\dw^2)=\dw^2\cha_s(\a_s\dw^{-2})$. Thus we 
have $\cha_s(\a_s(\dw^{-2}))=1$ that is $\cha_s(\a_s(w(t^q)t))=1$. Since $w(\a_s)=\a_s$ it follows that
$\cha_s(\a_s(t^{q+1}))=1$ hence $(\cha_s(-\a_s(t)))^{q+1}=1$. Thus (a) holds.

From our assumptions we have that $w(y')=y'$ and $s(y')=y'{}\i$ for any $y'\in T_s$; 
since $s(t)t\i\in T_s$, it follows that $w(s(t)t\i)=s(t)t\i$. Moreover we have $w(\ds^2)=\ds^2$. Hence for
$y\in T_s$ we have
$$sw(t^qy^q)ty(\ds\dw)^2=s(w(t^q)t\dw^2)sw(y^q)s(t)\i ty\ds^2=y^{-q}s(t)\i ty\ds^2.$$
This equals $1$ if and only if $y^{q-1}=s(t)\i t\ds^2$.  This proves the first sentence of (b).
The second sentence of (b) follows from (a).

We prove (c). For $m'\in N^{\ph'}$ and $(xU,\ph(x)U)\in O_{m'}^\ps$ fixed,
the variety $Z_{xU,\ph(x)U}$ in 2.3 (with $n=\ds$) can be identified with
$$Z'_{xU,\ph(x)U}=\{x'U\in X;(xU,x'U)\in O_{\ds m},(x'U,\ph(x)U)\in O_{\ds\i}\}.$$
(We use Lemma 1.4 and the equality $|sw|=|w|+1$.)
By Lemma 1.6, $Z'_{xU,\ph(x)U}$ is an affine line if $m'=\ds m\ds\i$, is a point if
$m'=\ds m y$ for some $y\in T_s$ and is empty otherwise. Hence 
$\sha(Z_{xU,\ph(x)U}^\ps)$ is $q$ if $m'=\ds m\ds\i$, is $1$ if
$m'=\ds m y$ for some $y\in T_s$ and is $0$ otherwise. 
Now (c) follows from (a),(b).

\proclaim{Lemma 2.7} Assume that $s\in S$, $w\in I$, $m\in N(w)$ and that
$sw=ws$, $|ws|<|w|$. Write $m=\dw t$ where $t\in T$ satisfies $w(t^q)t\dw^2=1$.

(a) For $y\in T_s$ we have $\ds m\ds\i y\in N(w)$ if and only if $y^{q-1}=1$.

(b) We have $s(t)t\i\e_s=m\i\ds m\ds\in T_s^\ph$.

(c) For $y\in T_s$ we have $\ds my\in N(sw)$ if and only if
$y^{q+1}=s(t)t\i\e_s=m\i\ds m\ds$. There are exactly $q+1$ such $y$; they are all automatically in $T_s^\Ph$.

(d) We have
$$k_{\ds}\th_m=q\sum_{y\in T_s;y^{q+1}=m\i\ds m\ds}\th{\ds my}+
\th_{\ds m\ds\i}+(q+1)\sum_{y\in T_s;y^{q-1}=1,y\ne1}\th_{\ds m\ds\i y}.$$
\endproclaim
We prove (a). We have
$$\align&\ph(\ds m\ds\i y)\ds m\ds\i y=\ds\ph(m)\ds\i y^q\ds m\ds\i y=\ds m\i y^{-q}m\ds\i y
\\&=\ds w(y^{-q})\ds\i y=\ds y^q\ds\i y=y^{-q}y=y^{1-q}.\endalign$$
This proves (a).

The equality in (b) is easily checked.
We have $s(t)t\i\e_s\in T_s$. To prove (b) it remains to show that $(s(t)t\i\e_s)^{q-1}=1$. 
We have $\ds\i\dw^2=\dw^2\ds\i$ hence 
$\dw^2=\ds\i\dw^2\ds=s(\dw^2)=\dw^2\cha_s(\a_s\dw^{-2})$. Thus we have
$\cha_s(\a_s(\dw^{-2}))=1$ that is $\cha_s(\a_s(w(t^q)t))=1$. Since $w(\a_s)=\a_s\i$ it follows that
$\cha_s(\a_s(t^{-q+1}))=1$ hence $(\cha_s(-\a_s(t)))^{-q+1}=1$. Thus (b) holds.

We prove (c). We have
$$\align&\ph(\ds my)\ds my=\ds\ph(m)y^q\ds my=\ds m\i y^q\ds my=\ds t\i\dw\i y^q\ds \dw ty\\&=
\ds t\i\dw\i y^q\dw\ds ty=\ds t\i w(y^q)\ds ty=\ds t\i y^{-q}\ds ty\\&=
s(t\i y^{-q})\e_sty=y^{q+1}s(t\i)t\e_s.\endalign$$ 
This proves the first sentence of (c). The second sentence of (c) follows from (b).

We prove (d). For $m'\in N^{\ph'}$ and $(xU,\ph(x)U)\in O_{m'}^\ps$ fixed, the variety $Z_{xU,\ph(x)U}$ in 2.3 
(with $n=\ds$) is 

(i) an affine line if $m'=\ds my$ for some $y\in T_s$ such that $\ds my\in N(sw)$, 

(ii) an affine line minus a point if $m'=\ds m\ds\i y$ with $y\in T_s-\{1\}$,

(iii) a union of two affine lines with one point in common if $m'=\ds m\ds\i$.
\nl
This is a geometric reinterpretation (and refinement) of the formula 1.7(e), in 
which the number of $\Ph$-fixed points on these
varieties enter; this number is $Q$ in case (i), is $Q-1$ in case (ii) and is $2Q-1$ in case (iii). 
It is enough to show that the number  of $\ps$-fixed points on $Z_{xU,\ph(x)U}$ is
$q$ in case (i), is $q+1$ in case (ii) and is $1$ in case (iii). This is verified directly by calculation in each
case. (In case (iii), $\ps$ interchanges the two lines, keeping fixed the point common to the two lines.)
We give the details of the calculation assuming that $G=SL_2(\kk)$, $T$ is the diagonal matrices, $TU$ is the 
upper triangular matrices, $\ds=\left(\sm 0&-1\\1&0\esm\right)$ and $\ph$ raises each matrix entry to the 
$q$-th power. We have $N^{\ph'}=\{M_a;a\in F_Q^*;a^q+a=0\}\sqc\{M'_a;a\in F_Q^*;a^{q+1}=1\}$ where 
$M_a=\left(\sm 0&-a\i\\a&0\esm\right)$, $M'_a=\left(\sm a&0\\0&a\i\esm\right)$. We must show: 

{\it if $x\in G$, $x\i\ph(x)=M_a$ then
$\sha(yU\in G/U; y\i\ph(y)\in UM_bU, x\i y\in U\ds U)=1+\d_{a,b}q$ (here $a^q+a=0,b^q+b=0$);}

{\it if $x\in G$, $x\i\ph(x)=M_{a'}$ then
$\sha(yU\in G/U; y\i\ph(y)\in UM_bU, x\i y\in U\ds U)=q$ (here $a'{}^{q+1}+a'=0$, $b^{q+1}=0$).}
\nl
Setting $y=xD$ we see that we must show that if $b^q+b=0$, then:

if $a^q+a=0$  then $\sha(DU\in(U\ds U)/U;D\i M_a\ph(D)\in UM_bU)=1+(1-\d_{a,b})q$;

if $a'{}^{q+1}=1$  then $\sha(DU\in (U\ds U)/U;D\i M_{a'}\ph(D)\in UM_bU)=q$.
\nl
Equivalently, we must show that if $b^q+b=0$, then:

(e) if $a^q+a=0$  then $\sha(d\in F_Q;d^{q+1}a-a\i=b)=1+(1-\d_{a,b})q$;

(f) if $a'{}^{q+1}=1$  then $\sha(d\in F_Q;-a'd^q+a'{}\i d=b)=q$.
\nl
If $a=b$, the equation in (e) is $d^{q+1}=0$ which has one solution namely $d=0$. If $a\ne b$ the
 equation in (e) is $d^{q+1}=ba\i+a^{-2}$. Here $(ba\i+a^{-2})^q=ba\i+a^{-2}\ne0$. Hence the equation in (e)
has exactly $q+1$ solutions. Setting $d'=a'{}\i d$, the equation in (f) is $-d'{}^q+d'=b$ and this has
exactly $d$ solutions in $F_Q$ since $b^q+b=0$. This completes the proof.

\subhead 2.8\endsubhead
Let $T(w)^*=\Hom(T(w),\CC^*)$. Since $e_w$ is surjective (see 2.1(b)), the map $T(w)^*@>>>\fs$, 
$\z\m\z e_w$ is an injective homomorphism. Let $\fs_w$ be the image of this homomorphism.
We have $\fs_w=\{\nu\in\fs;w(\nu)\nu^q=1\}$. 
Note that if $w\in W_2,z\in W$, then $z(\fs_w)=\fs_{zwz\i}$. 

For $\nu\in\fs_w$ we denote by $\unu_w$ the element of $T(w)^*$ such that $\nu=\unu_w e_w$.
We set $K_w=\ker(e_w)$.

For any $w\in I, n\in N(w)$ and $\nu\in\fs_w$ we define $a'_{n,\nu}\in\cf'$ by
$$a'_{n,\nu}=\sum_{t\in T(w)}\unu_w(t)\th_{nt}=|K_w|\i\sum_{\t\in T^\Ph}\nu(\t)\th_{\t n\t^{-q}}.$$
To verify the last equality we note that the sum over $t\in T(w)$ is equal to
$|K_w|\i\sum_{\t\in T^\Ph}\unu_w(e_w(\t))\th_{ne_w(\t)}$. We show:

(a) {\it If $w\in I,n\in N(w),\t\in T^\Ph,t\in T(w)$ and $\nu\in\fs_w$ then 
$a'_{nt,\nu}=\unu_w(t\i)a'_{n,\nu}$ and $a'_{\t n\t^{-q},\nu}=\nu(\t\i)a'_{n,\nu}$. 
In particular, the line spanned by $a'_{n,\nu}$ depends only on $w,\nu$ and not on $n$.}
\nl
Indeed, we have 
$$a'_{nt,\nu}=\sum_{t'\in T(w)}\unu_w(t')\th_{ntt'}=\sum_{t''\in T(w)}\unu_w(t''t\i)\th_{nt''}=\unu_w(t\i)a'_{n,\nu},$$
$$\align&a'_{\t n\t^{-q},\nu}\\&=|K_w|\i\sum_{\t'\in T^\Ph}\nu(\t')\th_{\t'\t n\t\t'{}^{-q}}
|K_w|\i\sum_{\t_1\in T^\Ph}\nu(\t_1\t\i)\th_{\t_1 n\t_1^{-q}}=\nu(\t\i)a'_{n,\nu}.\endalign$$
This proves (a).

From 2.1, 2.2(a),(b), we see that:

(b) {\it if $\{t_w;w\in W_2\}$ is a collection of elements in $T$ such that $\dw t_w\in N(w)$ for all $w\in W_2$ then
$\{a'_{\dw t_w,\nu};w\in W_2,\nu\in\fs_w\}$ is a $\CC$-basis of $\cf'$.}
\nl
For $\nu\in\fs,w\in I,n\in N(w),\nu'\in\fs_w$ we show:
$$1_\nu a'_{n,\nu'}=\d_{\nu,\nu'}a'_{n,\nu'}.\tag c$$
Indeed, we have
$$\align&1_\nu a'_{n,\nu'}=|K_w|\i|T^\Ph|\i\sum_{\t\in T^\Ph}\nu(\t)k_\t
\sum_{\t'\in T^\Ph}\nu'(\t')\th_{\t' n\t'{}^{-q}}\\&
=|K_w|\i|T^\Ph|\i\sum_{\t\in T^\Ph}\nu(\t)\sum_{\t'\in T^\Ph}\nu'(\t')\th_{\t\t' n\t'{}^{-q}t^{-q}}.\endalign$$
Setting $\t\t'=\t_1$ we obtain
$$\align&1_\nu a'_{n,\nu'}=|K_w|\i|T^\Ph|\i\sum_{\t_1\in T^\Ph}\nu'(\t_1)  
\sum_{\t\in T^\Ph}\nu(\t)\nu'(\t\i)\th_{\t_1n\t_1^{-1}{-q}}\\&
=\d_{\nu,\nu'}|K_w|\i\sum_{\t_1\in T^\Ph}\nu'(\t_1) \th_{\t_1n\t_1^{-1}{-q}}
=\d_{\nu,\nu'}a'_{n,\nu'}.\endalign$$
This proves (c).

For $s\in S$, $w\in W,n\in N(w),\nu\in\fs_w$, we have (using (c)):
$$\ct_sa'_{n,\nu}= q\i[\nu,\cha_s]k_{\ds}a_{n,\nu}.\tag d$$

\proclaim{Lemma 2.9}Let $s\in S,w\in I,n\in N(w),\nu\in\fs_w$. Note that $s\nu\in\fs_{sws}$. Assume that 
$sw\ne ws,|sw|>|w|$. We have
$$\ct_sa'_{n,\nu}=q\i[\nu,\cha_s]a'_{\ds n\ds\i,s\nu}.$$
\endproclaim
Using 2.8(d) we see that it is enough to show
$$k_{\ds}a'_{n,\nu}=a'_{\ds n\ds\i,s\nu}.$$
Using 2.5 and the equality $|K_w|=|K_{sws}|$ we see that
$$\align&k_{\ds}a'_{n,\nu}=|K_w|\i\sum_{\t\in T^\Ph}\nu(\t)k_{\ds}\th_{\t n\t^{-q}}\\&
=|K_w|\i\sum_{\t\in T^\Ph}\nu(\t)\th_{\ds\t n\t^{-q}\ds\i}\\&=
|K_{sws}|\i\sum_{\t'\in T^\Ph}\nu(s(\t'))\th_{\t'\ds n\ds\i\t'{}^{-q}}=a'_{\ds n\ds\i,s\nu}.\endalign$$
The lemma is proved.

\proclaim{Lemma 2.10}Let $s\in S,w\in I,n\in N(w),\nu\in\fs_w$. Assume that $sw=ws$, $|sw|>|w|$.
If $s\in W_\nu$ we set $\D=1$; if $s\n W_\nu$ we set $\D=0$.
Note that we have $s\nu\in\fs_w$; moreover, if $\D=1$ then $s\nu=\nu\in\fs_{sw}$. 
We set $z=\ds n\i\ds n\in T_s$, see 2.6(a). We have $z^{q+1}=1$, see 2.6(a). We have
$$\ct_sa'_{n,\nu}=a'_{\ds n\ds\i,s\nu}\text{ if }\D=0,$$
$$\ct_sa'_{n,\nu}=a'_{n,\nu}+(q\i+1)a'_{\ds nu,\nu}\text{ if }\D=1$$
where $u\in T_s^\Ph$ is such that $u^{q-1}=z$ (see 2.6(b)).
\endproclaim
Using 2.6(c) and 2.8(d) we have $\ct_sa'_{n,\nu}=A+B$ where
$$A=|K_w|\i\sum_{\t\in T^\Ph}\nu(\t)\th_{\ds\t nt^{-q}\ds\i}.$$
$$B=q\i|K_w|\i\sum_{\t\in T^\Ph,y\in T_s;y^{q-1}=\ds\t^qn\i\t\i\ds\t n\t^{-q}}\nu(\t)\th_{\ds\t n\t^{-q}y}.$$
We have used that $\nu(\e_s)=1$ (hence $[\nu,\cha_s]=1$). Indeed, we have $\nu(\e_s)=\unu_w(e_w(\e_s))=
\unu_w(w(\e_s)\e_s)=\unu_w(1)=1$ since $w(\e_s)=\e_s$.

In the sum $A$ we set $\t'=s(\t)$. We get
$$A=|K_w|\i \sum_{\t'\in T^\Ph}(s\nu)(\t')\th_{t'\ds n\ds\i t'{}^{-q}}=a'_{\ds n\ds\i,s\nu}.$$
We now show that if $\D=1$ then 
$$a'_{\ds n\ds\i,s\nu}=a'_{n,\nu}.$$
We write $n=\dw t$ with $t\in T$.
We have $\ds n\ds\i=\ds\dw t\ds\i=\dw\ds t\ds\i=nt\i s(t)$.
By 2.6(a) we have $(t\i s(t))^{q+1}=1$. Since $t\i s(t)\in T_s$ we have $t_1\i s(t)=t_1^{q-1}$ with $t_1\in T_s^\Ph$.
Thus we have $\ds n\ds\i=nt_1^{q-1}$ hence $a'_{\ds n\ds\i,\nu}=a'_{nt_1^{q-1},\nu}=
a'_{t_1\i nt_1^q,\nu}=a'_{n,\nu}$ since $\nu(t_1)=1$. This proves our claim.

We now consider the sum $B$. In that sum we have
$$\align&\ds\t^qn\i\t\i\ds\t n\t^{-q}=s(\t^q)\ds n\i\ds s(\t)\i\t n\t^{-q}\\&=
\ds n\i\ds ns(\t)\i\t\t^{-q}s(\t^q)=z(\t s(\t)\i)^{1-q}.\endalign$$
Thus we have
$$B=q\i|K_w|\i\sum_{(\t,y)\in \cy}\nu(\t)\th_{\ds nw(\t)\t^{-q}y}$$
where $\cy=\{(\t,y)\in T^\Ph\T T_s;y^{q-1}=z(\t s(\t)\i)^{1-q}\}$. 
Let $\cy'=\{(\t',u)\in T^\Ph\T(T_s^\Ph);u^{q-1}=z\}$.
The map $\x:\cy'@>>>\cy$, $(\t',u)\m(s(\t'),s(\t')^q\t'{}^{-q}u)$ is a well defined bijection.
Now the sum $B$ can be written in terms of this bijection as follows:
$$B=q\i|K_w|\i\sum_{(\t',u)\in \cy'}\nu(s(\t'))\th_{\ds nw(s(\t'))\t'{}^{-q}u}.$$
We have a free action of $T_s^\Ph$ on $\cy'$ given by $e:(\t',u)\m(\t's(e),ue^{-q-1})$. 
Note that the quantity $\th_{\ds nw(s(\t'))\t'{}^{-q}u)}$ is constant on the orbits of this action.
Hence if $\cy'_0$ is a set of representatives for the $T_s^\Ph$-orbits on $\cy'$ we have
$$B=q\i|K_w|\i\sum_{(\t',y)\in\cy'_0,e\in T_s^\Ph}(s\nu)(\t')\nu(e)\th_{\t'\ds n\t'{}^{-q}u}.$$
Note that $\sum_{e\in T_s^\Ph}\nu(e)=\d(q^2-1)$. In particular, if $\D=0$ we have $B=0$.
We now assume that $\D=1$.
For any $u\in T_s^\Ph$ such that $u^{q-1}=z$ we set
$$B_u=q\i|K_w|\i\sum_{\t'\in T^\Ph}(s\nu)(\t')\th_{\t'\ds n\t'{}^{-q}u}.$$
We have $B=\sum_{u\in T_s^\Ph;u^{q-1}=z}B_u$.
For any $u$ as above and any $e\in T_s^\Ph$ we have $B'_{ue^{-1-q}}=B'_u$ since $\nu(e)=1$.
If $u,u'$ in  $T_s^\Ph$ are such that $u^{q-1}=u'{}^{q-1}=z$, we have $u'=u\ti e$ where
$\ti e\in T_s^\Ph$ satisfies $\ti e^{q-1}=1$. Hence we have $\ti e=e^{-q-1}$ for some
$e\in T_s^\Ph$ so that $u'=ue^{-q-1}$. Thus we have $B_{u'}=B_u$. We see that
$B=(q-1)B_u$ where $u\in T_s^\Ph$ is such that $u^{q-1}=z$. 
We have $B_u=q\i|K_{sw}||K_w|\i a'_{\ds nu,\nu}$.
It remains to show that $(q-1)|K_{sw}||K_w|\i=q+1$ or equivalently, that 
$|T(sw)||T(w)|\i=(q-1)(q+1)\i$. This follows from the following fact: there exists $c,c'$ in $\NN$ such that
$|T(w)|=(q-1)^c(q+1)^{c'}$, $|T(sw)|=(q-1)^{c+1}(q+1)^{c'-1}$. The lemma is proved.

\proclaim{Lemma 2.11}Let $s\in S,w\in I,n\in N(w),\nu\in\fs_w$. Assume that $sw=ws$, $|sw|<|w|$,
$s\n W_\nu$. Note that $s\nu\in\fs_w$. We have
$$\ct_sa'_{n,\nu}=-a'_{\ds n\ds\i,s\nu}.$$
\endproclaim
Using 2.7(d) and 2.8(d) we have $\ct_sa'_{n,\nu}=A+B$ where 
$$A=q\i|K_w|\i\sum_{\t\in T^\Ph,y\in T_s;y^{q-1}=1}c_y\nu(\t)\th_{\ds\t n\t^{-q}\ds\i y},$$
where $c_y=q+1$ if $y\ne1$, $c_y=1$ if $y=1$ and
$$B=|K_w|\i \sum_{\t\in T^\Ph,y\in T_s;y^{q+1}=\t^q n\i\t\i\ds \t n\t^{-q}\ds}\nu(\t)
\th_{\ds\t n\t^{-q}y}.$$
We have used that, as in the proof of 2.10, we have $\nu(\e_s)=1$ (hence $[\nu,\cha_s]=1$).
In the sum $A$ we set $\t'=s(\t)$. We get
$$A=q\i|K_w|\i\sum_{\t'\in T^\Ph,y\in T_s;y^{q-1}=1}c_y(s\nu)(\t')\th_{\t'\ds n\ds\i\t'{}^{-q}y}.$$
Fr $y\in T_s$ such that $y^{q-1}=1$ we can find $y'\in T_s$ such that $y'{}^{q+1}=y$ (there are $q+1$ such $y'$)
and we have automatically $y'\in T^\Ph$. Thus we have
$$\align&A=q\i|K_w|\i(q+1)\i\sum_{\t'\in T^\Ph,y'\in T_s^\Ph}c_{y'{}^{-q-1}}(s\nu)(\t')
\th_{\t'\ds n\ds\i\t'{}^{-q}y'{}^{-q-1}}\\&
=q\i|K_w|\i(q+1)\i\sum_{\t'\in T^\Ph,y'\in T_s^\Ph}c_{y'{}^{-q-1}}(s\nu)(\t')
\th_{y'\t'\ds n\ds\i y'{}^{-q}\t'{}^{-q}}.\endalign$$
With the change of variable $\t'y'=\t''$ we obtain
$$A=q\i|K_w|\i(q+1)\i\sum_{\t''\in T^\Ph,y'\in T_s^\Ph}c_{y'{}^{-q-1}}(s\nu)(\t'')\nu(y')
\th_{\t''\ds n\ds\i\t''{}^{-q}}.$$
(We have used that $s(y')=y'{}\i$.) Using our assumption that $s\n W_\nu$, we have
$$\align&\sum_{y'\in T_s^\Ph}c_{y'{}^{-q-1}}\nu(y')\\&
=\sum_{y'\in T_s^\Ph;y'{}^{q+1}=1}\nu(y')+(q+1)\sum_{y'\in T_s^\Ph;y'{}^{q+1}\ne1}\nu(y')\\&
=(q+1)\sum_{y'\in T_s^\Ph}\nu(y')-q\sum_{y'\in T_s^\Ph;y'{}^{q+1}=1}\nu(y')\\&
=-q\sum_{y'\in T_s^\Ph;y'{}^{q+1}=1}\nu(y')=-q\sum_{y'\in T_s^\Ph;y'{}^{q+1}=1}\unu_w(y'{}^{-q-1})\\&
=-q\sha(y'\in T_s^\Ph;y'{}^{q+1}=1)=-q(q+1).\endalign$$
It follows that
$$A=q\i|K_w|\i(q+1)\i(-q)(q+1)\sum_{\t''\in T^\Ph}\nu(\t'')\th_{\t''\ds n\ds\i\t''{}^{-q}}
=-a'_{\ds n\ds\i,s\nu}.$$
It remains to prove that $B=0$. We set $z=n\i\ds n\ds\in T_s$, see 2.7(b). In the sum $B$ we have
$$\align&\t^q n\i\t\i\ds\t n\t^{-q}\ds=\t^q n\i\t\i s(\t)      \ds n\ds s(\t)^{-q}\\&=
 \t^q\t s(\t)\i n\i   \ds n\ds s(\t)^{-q}=z\t^q\t s(\t)\i s(\t)^{-q}=z(\t s(\t)\i)^{q+1}.\endalign$$
Thus we have
$$B=q\i|K_w|\i q\sum_{(\t,y)\in\cz}\nu(\t)\th_{\ds\t n\t^{-q}y}$$
where $\cz=\{(\t,y)\in T^\Ph\T T_s;y^{q+1}=z(\t s(\t)\i)^{q+1}\}$.
The group $T_s^\Ph$ acts freely on $\cz$ by $e:(\t,y)\m(\t e,ye^{q+1})$. 
(We must show that the equation $y^{q+1}=z(\t s(\t)\i)^{q+1}$ implies
$(ye^{q+1})^{q+1}=z(\t es(\t e)\i)^{q+1}$; it is enough to show that $e^{(q+1)^2}=e^{2(q+1)}$ and this follows
from $e^{q^2-1}=1$.)
We show that the last sum restricted to any $T_s^\Ph$-orbit is zero. Since 
$\th_{\ds\t n\t^{-q}y}$ is constant on any  $T_s^\Ph$-orbit
it is enough to show that $\sum_{e\in T_s^\Ph}\nu(e)=0$; this follows from our assumption that $s\n W_\nu$.
We deduce that $B=0$. The lemma is proved.

\subhead 2.12\endsubhead
For $w\in I$ let $||w||$ be the dimension of the $-1$ eigenspace of the linear map induced by $w$ on the
real vector space $\RR\ot Y$. We have $|w|=||w||\mod2$. For $w\in N(w),\nu\in\fs_w$ we set
$$\ta_{n,\nu}=q^{-(|w|+||w||)/2}a'_{n,\nu}\in\cf'.$$
We have the following result.

\proclaim{Lemma 2.13}Let $s\in S,w\in W_2,n\in N(w),\nu\in\fs_w$. Write $n=\dw t$ where $t\in T$. We have:

(a) $\ct_s\ta_{n,\nu}=[\nu,\cha_s]\ta_{\ds n\ds\i,s\nu}$ if $sw\ne ws,|sw|>|w|$;

(b) $\ct_s\ta_{n,\nu}=\ta_{\ds n\ds\i,s\nu}$ if $sw=ws, |sw|>|w|,s\n W_\nu$;

(c) $\ct_s\ta_{n,\nu}=\ta_{n,\nu}+(q+1)\ta_{\ds nu,\nu}$
(where $u\in T_s^\Ph$ is such that $u^{q-1}=\ds n\i\ds n=s(t)\i t\e_s$, see 2.6(a),(b)) if $sw=ws$, $|sw|>|w|,s\in W_\nu$;

(d) $\ct_s\ta_{n,\nu}=-\ta_{\ds n\ds\i,s\nu}$ if $sw=ws$, $|sw|<|w|$, $s\n W_\nu$;
\endproclaim
(a) is a reformulation of 2.9; (b),(c) are reformulations of 2.10; (d) is a reformulation of 2.11. 

\head 3. Proof of Theorem 0.4\endhead
\subhead 3.1\endsubhead
We preserve the setup of 1.1.
Let $L$ be the subgroup of $Y$ generated by $\{\cha_s;s\in S\}$. 
Let $S'$ be a {\it halving} of $S$ that is a subset $S'$ of $S$ such that $s_1s_2=s_2s_1$ whenever $s_1,s_2$ in $S$ are 
both in $S'$ or both in $S-S'$. (Such $S'$ always exists.) 
Let $W_2@>>>Y$, $w\m r_w$ and $W_2@>>>L/2L$, $w\m b_w=b_w^{S'}$ be the maps defined in \cite{\LIF, 0.2, 0.3}. From {\it loc. cit.}
and from the proof of \cite{\LIF, 1.14(a)} we have:

(i) $r_1=0$, $r_s=\cha_s$ for any $s\in S$, $b_1=0$, $b_s=\cha_s$ for any $s\in S'$, $b_s=0$ for any $s\in S-S'$;

(ii) for any $w\in W_2,s\in S$ such that $sw\ne ws$ we have $s(r_w)=r_{sws}$, $s(b_w)=b_{sws}+\cha_s$;

(iii) for any $w\in W_2,s\in S$ such that $sw=ws$ we have $r_{sw}=r_w+\cn\cha_s$,
$b_{sw}=b_w+l\cha_s$ where $l\in\{0,1\}$, $\cn\in\{-1,0,1\}$.

(iv) for any $w\in W_2,s\in S$ such that $sw=ws,|sw|>|w|$ we have $s(r_w)=r_w$;

(v) for any $w\in W_2,s\in S$ such that $sw=ws$ we have $s(b_w)=b_w+(1-\cn)\cha_s$ where $\cn$ is as in (iii).

Moreover, by \cite{\LIF, 0.5}, 

(vi) if $c\in F_Q,c^{q-1}=\e$, the element $n_{w,c}=\dw r_w(c)b_w(\e)\in\k\i(w)$ belongs to $N(w)$.
\nl
Here $r_w(c)\in T$, $b_w(\e)\in T$ are obtained by evaluating a homomorphism $\kk^*@>>>Y$ at $c$ or $\e$. Note that $b_w(\e)=b_w(\e)\i$.
From \cite{\LIF, 1.18} we deduce:

(vii) in the setup of (iii) we have $\cn=(w:s)$.
\nl
The following equality complements (iv):

(viii) for any $w\in W_2,s\in S$ such that $sw=ws,|sw|<|w|$ we have $s(r_w)=r_w+2(w:s)\cha_s$.
\nl
Indeed, using (iii),(iv),(vii) we have

$s(r_w)=s(r_{sw}-(w:s)\cha_s)=r_{sw}+(w:s)\cha_s=r_w+2(w:s)\cha_s$.

\mpb

For any $w\in W_2$, any $c\in F_Q$ such that $c^{q-1}=\e$ and any $\nu\in\fs_w$ we set
$$a_{w,c,\nu}=\ta_{n_{w,c},\nu}.$$
This is well defined by (vi). By 2.8(b), 

(a) {\it for any $c$ as above, $\{a_{w,c,\nu};w\in W_2,\nu\in\fs_w\}$ is a $\CC$-basis of $\cf'$.}
\nl
In the remainder of this section we assume that 0.3(a) holfds. We have the following result.

\proclaim{Proposition 3.2}Let $s\in S,w\in W_2,\nu\in\fs_w$. Let $c$ be as in 3.1(vi).
We have 

(a) $\ct_sa_{w,c,\nu}=a_{sws,c,s\nu}$ if $sw\ne ws,|sw|>|w|,s\in W_\nu$;  

(b) $\ct_sa_{w,c,\nu}=a_{sws,c,s\nu}+(q-q\i)a_{w,c,\nu}$ if $sw\ne ws,|sw|<|w|,s\in W_\nu$;

(c) $\ct_sa_{w,c,\nu}=a_{w,c,\nu}+(q+1)a_{sw,c,\nu}$ if $sw=ws,|sw|>|w|,s\in W_\nu$;

(d) $\ct_sa_{w,c,\nu}=(1-q\i)a_{sw,c,\nu}+(q-q\i-1)a_{w,c,\nu}$ if $sw=ws,|sw|<|w|,s\in W_\nu$;

(e)  $\ct_sa_{w,c,\nu}=[\nu,\cha_s]a_{sws,c,s\nu}$ if $sw\ne ws,|sw|>|w|,s\n W_\nu$;

(f) $\ct_sa_{w,c,\nu}=[\nu,\cha_s]\i a_{sws,c,s\nu}$ if $sw\ne ws,|sw|<|w|,s\n W_\nu$;

(g) $\ct_sa_{w,c,\nu}=\un{s\nu}_w(\e_s^{1-(w:s)})a_{w,c,s\nu}$ if $sw=ws, |sw|>|w|,s\n W_\nu$;

(h) $\ct_sa_{w,c,\nu}=-\un{s\nu}_w(\e_s^{1-(w:s)})\un{s\nu}_w(c_s^{-2(w:s)})a_{w,c,s\nu}$ 
if $sw=ws,|sw|<|w|,s\n W_\nu$.
\endproclaim
This will be deduced in 3.3-3.8 from 2.13 with $n=n_{w,c}$ as in 3.1(vi), using the equality 
$\ta_{n't,\nu'}=\unu'_{w'}(t\i)\ta_{n',\nu'}$ where $w'\in W_2,n'\in N(w),\nu'\in\fs_{w'}, t\in T(w')$, which
follows from 2.8(a). 

\subhead 3.3\endsubhead
Assume that we are in the setup of 3.2(a) or 3.2(e). Using  2.13(a) and 3.1(ii) we obtain
$$\align&\ct_sa_{w,c,\nu}=[\nu,\cha_s]\ta_{\ds \dw r_w(c)b_w(\e)  \ds\i,s\nu}\\&=
[\nu,\cha_s]\ta_{\ds \dw\ds \ds\i r_w(c)b_w(\e)  \ds\i,s\nu}\\&=
[\nu,\cha_s]\ta_{n_{sws,c} r_{sws}(c)\i b_{sws}(\e)\i  ds\i r_w(c)b_w(\e)  \ds\i,s\nu}\\&=
[\nu,\cha_s]\un{s\nu}_{sws}(\ds b_w(\e)r_w(c)\i \ds r_{sws}(c)b_{sws}(\e))a_{sws,c,s\nu}\\&
=[\nu,\cha_s]\un{s\nu}_{sws}(b_{sws}(\e)\e_sr_{sws}(c)\i\e_sr_{sws}(c)b_{sws}(\e))a_{sws,c,s\nu}
=[\nu,\cha_s]a_{sws,c,s\nu}.\endalign$$
This proves 3.2(e). 
Now 3.2 follows also since in that case we have $[\nu,\cha_s]=1$.
(It is enough to show that $\nu(\e_s)=1$. This follows from $s\in W_\nu$.) This proves 3.2(a).

\subhead 3.4\endsubhead
Assume that we are in the setup of 3.2(g). Using 2.13(b), 3.1(iv),(v),(vii), we obtain
$$\align&\ct_sa_{w,c,\nu}=\ta_{\ds\dw r_w(c)b_w(\e)\ds\i,s\nu}=\ta_{\dw s(r_w(c)b_w(\e)),s\nu}\\&=
\ta_{\dw r_w(c)b_w(\e) \e_s^{1-(w:s)},s\nu}=\un{s\nu}_w(\e_s^{1-(w:s)})a_{w,c,s\nu}.\endalign$$
This proves 3.2(g).

\subhead 3.5\endsubhead
Assume that we are in the setup of 2.13(c) with $n=n_{w,c}$. Using 3.1(iv), (v), (vii), we have
$$\align&u^{q-1}=s(r_w(c)b_w(e))\i r_w(c)b_w(e)\e_s=s(b_w(e))\i b_w(e)\e_s\\&=\e_s^{1-(w:s)}\e_s=\e_s^{(w:s)}.
\tag a\endalign$$
For $l\in\{0,1\}$ we show:
$$\unu_{sw}(c_s^{(w:s)}\e_s^lu\i)=1.\tag b$$
Since $\nu$ is $1$ on $T_s^\Ph$, $\unu_{sw}$ must be trivial on 
$e_{sw}(T_s^\Ph)$ that is on the image of $T_s^\Ph@>>>T_s^\Ph$, $t\m t^{q+1}$ which is the
same as $\{t'\in T_s;t'{}^{q-1}=1\}$. Since 
$c_s^{(w:s)}\e_s^lu\i\in T_s$, it is enough to show that
$$(c_s^{(w:s)}\e_s^lu\i)^{q-1}=1.\tag c$$
Using (a) and the equations $c^{q-1}=\e$, $\e^{q-1}=1$, we see that the left hand side of (c) is
$\e_s^{(w:s)}\e_s^{-(w:s)}=1$. This completes the proof of (b).

We now assume that we are in the setup of 3.2(c) (which is the same as the setup of 2.13(c) with $n=n_{w,c}$).
From 2.13(c) we deduce using (b) and 3.1(iii) that for some $l\in\{0,1\}$ we have
$$\align&\ct_sa_{w,c,\nu}-a_{w,c,\nu}=(q+1)\ta_{\ds\dw r_w(c)b_w(\e)u,\nu}\\&=
(q+1)\ta_{\ds\dw r_{sw}(c)b_{sw}(\e) r_{sw}(c)\i b_{sw}(\e)r_w(c)b_w(\e)u,\nu}\\&=
(q+1)\unu_{sw} (r_{sw}(c)b_{sw}(\e) r_w(c)\i b_w(\e)u\i)a_{sw,c,\nu}\\&
=(q+1)\unu_{sw} (c_s^{(w:s)}\e_s^lu\i)a_{sw,c,\nu}=(q+1)a_{sw,c,\nu}.\endalign$$
This completes the proof of 3.2(c).

\subhead 3.6\endsubhead
Assume that we are in the setup of 3.2(h). From 2.13(d) we deduce using 3.1(viii):
$$\align&\ct_sa_{w,c,\nu}=-\ta_{\ds \dw r_w(c)b_w(\e)\ds\i,s\nu}
=-\ta_{\dw s(r_w(c)b_w(\e)),s\nu}\\&=
-\ta_{\dw r_w(c)b_w(\e)c_s^{2(w:s)}\e_s^{1-(w:s)},s\nu}=\un{s\nu}_w(c_s^{-2(w:s)}\e_s^{1-(w:s)})a_{w,c,s\nu}\\&=
\un{s\nu}_w(\e_s^{1-(w:s)}\un{s\nu}_w(c_s^{-2(w:s)})a_{w,c,s\nu}.\endalign$$
This proves 3.2(h).

\subhead 3.7\endsubhead    
Assume that $sw\ne ws$, $|sw|<|w|$. Then 3.2(a),(e) are applicable with $sws,s\nu$ instead of $w,\nu$
so that
$$\ct_sa_{sws,c,s\nu}=[s\nu,\cha_s]a_{w,c,\nu}.$$
We apply $\ct_s\i$ to both sides; we obtain
$$\ct_s\i a_{w,c,\nu}=\ct_s\i1_\nu a_{w,c,\nu}=[s\nu,\cha_s]\i a_{sws,c,s\nu}.$$
Using 1.8(i) we deduce
$$\ct_sa_{w,c,\nu}-\d(q-q\i)a_{w,c,\nu}=[s\nu,\cha_s]\i a_{sws,c,s\nu}$$
where $\D=1$ if $s\in W_\nu$, $\D=0$ if $s\n W_\nu$. Thus,
This proves 3.2(b),(f). (We use that $[s\nu,\cha_s]=[\nu,\cha_s]$ is $1$ when $s\in W_\nu$ since $\nu(\e_s)=1$ in that case.)

\subhead 3.8\endsubhead    
Assume that $s,w,\nu$ are as in 3.2(d). Then 3.2(c) is applicable to $sw,\nu$ instead of $w,\nu$ and gives:
$$\ct_sa_{sw,c,\nu}=a_{sw,c,\nu}+(q+1)a_{w,c,\nu}.\tag a$$
We apply $\ct_s$ to (a). We obtain
$$\ct_s\ct_sa_{sw,c,\nu}=\ct_sa_{sw,c,\nu}+(q+1)\ct_sa_{w,c,\nu}.$$
Using 1.8(h) we deduce
$$a_{sw,c,\nu}+(q-q\i)\ct_sa_{sw,c,\nu}=\ct_sa_{sw,c,\nu}+(q+1)\ct_sa_{w,c,\nu}$$
hence, using (a):
$$a_{sw,c,\nu}+(q-q\i-1)(a_{sw,c,\nu}+(q+1)a_{w,c,\nu})=(q+1)\ct_sa_{w,c,\nu}.$$
Dividing by $q+1$ we get 3.2(d). This completes the proof of Proposition 3.2.

\subhead 3.9\endsubhead    
We choose a generator $\g$ of the cyclic group $F_Q^*$ so that we have an isomorphism
$$\ZZ/(Q-1)\ZZ@>\si>>F_Q^*\tag a$$
 which takes $1$ to $\g$. 

Let $z\in\ZZ$ be as in Theorem 0.2. Let $c=\g^{z(q+1)/2}\in F_Q^*$.
(If $p=2$ so that $(q+1)/2$ is not an integer, this is interpreted as a square root of $\g^{z(q+1)}$
which is uniquely defined.) If $p\ne2$ we have $c^{q-1}=\g^{z(q^2-1)/2}=\e$ by the choice of $z$.
If $p=2$ then $(c^{q-1})^2=(c^2)^{q-1}=\g^{z(q^2-1)}=1$ hence $c^{q-1}=1=\e$. Thus in any case we
have $c^{q-1}=\e$.

We have an isomorphism of groups
$F_Q^*\ot Y@>\si>>T^\Ph$, $z\ot y\m y(z)$. Using (a) this can be viewed as an isomorphism
of groups $(\ZZ/(Q-1)\ZZ)\ot Y@>\si>>T^\Ph$; it takes $n\ot y$ to $y(\g^n)$. We have a pairing 
$$(,):((\ZZ/(Q-1)\ZZ)\ot Y)\T \bX_q@>>>\CC^*$$
given by 
$$(d\ot y,\fra{a}{Q-1}\ot x)=\exp(2\pi\sqrt{-1}\fra{da}{Q-1}\lf y,x\rf)$$
where $y\in Y,x\in X$, $a\in\ZZ$, $d\in\ZZ$. This pairing identifies
$\bX_q$ with \lb
$\Hom((\ZZ/(Q-1)\ZZ)\ot Y,\CC^*)=\Hom(T^\Ph,\CC^*)=\fs$.
This identification is compatible with the natural $W$-actions on $\bX_q$ and $\fs$; it induces
an identification $\tX_q=\{(w,\nu);w\in W_2,\nu\in\fs_w\}$. 
Thus, the basis 3.1(a) of $\cf'$ can be naturally indexed by the elements of $\tX_q$.
We shall interpret the quantities
$$[\nu,\cha_s],\un{s\nu}_w(\e_s^{1-(w:s)}),\un{s\nu}_w(c_s^{-2(w:s)})$$
 which appear in 3.2 in terms of the corresponding parameter in $\tX_q$.
Assume that $(w,\nu)\in W_2\T\fs$ (with $\nu\in\fs_w$) corresponds to $(w,\l)\in\tX_q$.
Then for any $s\in S$ we have
$$\nu(\cha_s(\g))=\exp(2\pi\sqrt{-1}\lf\cha_s,\l\rf).\tag b$$
We show:

(c) If $sw=ws,|sw|<|w|,s\n W_\nu$ then 
$$\un{s\nu}_w(c_s^{-2(w:s)})=\exp(2\pi\sqrt{-1}(w:s)z\lf\cha_s,\l\rf).$$
Let $\tc=\g^z$. We have $\tc_s^{q+1}=c_s^2$ hence
$$\un{s\nu}_w(c_s^{-2(w:s)})=\un{s\nu}_w((\tc_s^{-(w:s)})^{q+1})=\un{s\nu}_w(e_w(\tc_s^{-(w:s)}))=
(s\nu)(\tc_s^{-(w:s)})=\nu(\tc_s^{(w:s)}).$$
It remains to show:
$$\nu(\cha_s(\g^z)=\exp(2\pi\sqrt{-1}z\lf\cha_s,\l\rf).$$
This clearly follows from (b).

\mpb  

We show:

(d) If $sw=ws$ then $\un{s\nu}_w(\e_s^{1-(w:s)})=\d_{w,s\l,s}$.
\nl
If $p=2$, both sides are $1$. Thus we can assume that $p\ne2$. We must show that
$$\un{s\nu}_w(\e_s^{1-(w:s)})=\exp(2\pi\sqrt{-1}((q-e)/2)(1-(w:s))\lf\cha_s,s\l\rf)$$
where $e=|w|-|sw|=\pm1$. It is enough to show that
$$\un{s\nu}_w(\e_s)=\exp(2\pi\sqrt{-1}((q-e)/2)\lf\cha_s,s\l\rf)$$
We have $\e_s=(\g_s^{(q-e)/2})^{q+e}=e_w(\g_s^{(q-e)/2})$ so that
$$\un{s\nu}_w(\e_s)=\un{s\nu}_w(e_w(\g_s^{(q-e)/2}))=(s\nu)(\g_s^{(q-e)/2}).$$
Thus it is enough to show that
$$(s\nu)(\cha_s(\g))=\exp(2\pi\sqrt{-1}\lf\cha_s,s\l\rf)$$
This clearly follows from (b).

\mpb

We show:

(e) If $s\in S$ then $[\l,s]=[\nu,\cha_s]$.
\nl
If $p=2$ both sides are $1$. Thus we can assume that $p\ne2$. We must show that
we have $[\l,s]=1$ if and only if $[\nu,\cha_s]=1$ or that
$\exp(2\pi\sqrt{-1}(1/2)(Q-1)\lf\cha_s,s\l\rf)=1$ if and only if $\nu(\cha_s(\e))=1$
or (using (b)) that $\nu(\cha_s(\g))^{(1/2)(Q-1)}=1$ if and only if $\nu(\cha_s(\e))=1$.
This follows from the equality $\g^{(1/2)(Q-1)}=\e$.

From (b) and the definitions we see that:

(f) If $s\in S$ then we have $s\in W_\l$ if and only if $s\in W_\nu$.
\nl
We now see that Proposition 3.2 implies the truth of Theorem 0.4 in the special case
where $\kk$ is as in 1.1. But then Theorem 0.4 follows immediately for any $\kk$ as
in 0.1 such that the characteristic of $\kk$ is $0$ or $p$. This completes the proof
of Theorem 0.4.

\head 4. The generic case\endhead
\subhead 4.1\endsubhead
In this section we assume that $\kk=\CC$ and that 0.3(a) holds.
We have $\bX_1=\bX$. Hence $\tX_1=\{(w,\l)\in W_2\T\bX; w(\l)=-\l\}$.

Until the end of 4.2, we fix a $W$-orbit $\co$ in $\bX$ which is contained in the image of $X_\QQ$ under $X_K@>>>\bX$.
We can find an integer $\fe\ge1$ such that $\fe\lf y,\l\rf=0$ for any $y\in Y$ and any
$\l\in\co$.
We can write $\fe=\prod_{p\in\fP}p^{c_p}$ where $\fP$ is a finite set of prime numbers and $c_p\ge1$ are integers.
Let $\fP'$ be the set of prime numbers which do not divide $2\fe$. Note that $\fP\cap\fP'=\emp$. Hence
if $p\in\fP,p'\in\fP'$ then $p'$ is a unit in the ring $\ZZ/p^{c_p}\ZZ$ hence for some integer $a_p\ge1$
independent of $p'$ we have $p'{}^{a_p}=1$ in $\ZZ/p^{c_p}\ZZ$ that is $p^{c_p}$ divides $p'{}^{a_p}-1$.
Let $\cs$ be the set of all integers $z\ge1$ such that $z$ is divisible by $\prod_{\p\in\fP}a_p$.
Then for any $p\in\fP,p'\in\fP'$ and any $z\in\cs$, $p^{c_p}$ divides $p'{}^z-1$.
Hence for any $p'\in\fP'$ and any $z\in\cs$, $\fe$ divides $p'{}^z-1$.
Let $\fQ$ be the set of all numbers of the form $p'{}^z$ with $p'\in\fP',z\in\cs$. Then
we have $(q-1)\lf y,\l\rf=0$ for any $q\in\fQ$, any $y\in Y$ and any $\l\in\co$. Hence

(a) $(q-1)\l=0$ for any $q\in\fQ$ and any  $\l\in\co$. 
\nl
It follows that

(b) if $(w,\l)\in\tX_1$ and $\l\in\co$ then $(w,\l)\in\tX_q$ for any $q\in\fQ$.
\nl
Indeed, we have $w(\l)=-\l$ and we must show that $w(\l)=-q\l$. It is enough to show that
$q\l=\l$ and this follows from (a).

\subhead 4.2\endsubhead
Let $\tfQ$ be the set of squares of the numbers in $\fQ$. We have $\tfQ\sub\fQ$.
We now fix $q\in\tfQ$. We have $q=q'{}^2$ with $q'\in\fQ$. Note that $q=4t+1$ for some $t\in\NN$.
Let $(w,\l)\in\tX_1$ with $\l\in\co$ (so that $(w,\l)\in\tX_{q'}$ and $(w,\l)\in\tX_q$ by 4.1(b)) and let $s\in S$.
We show:

(a) $[\l,s]$ defined as in 0.2 in terms of $q$ is equal to $1$. 
\nl
Since $(w,\l)\in\tX_{q'}$ we have $\lf\cha_s,\l\rf=e'/(q'{}^2-1)$ with $e'\in\ZZ$.
Hence $\lf\cha_s,\l\rf=e/(q^2-1)$ with $e=e'(q'{}^2+1)$. Since $e$ is even we see that (a) holds.

\mpb

We show:

(b) If $sw=ws$, $|sw|>|w|$ then $\d_{w,\l;s}$ defined as in 0.3 in terms of $q$ is equal to $\d'_{w,\l;s}$ defined as in 0.5.
\nl
It is enough to show that $\exp(2\pi\sqrt{-1}((q+1)/2)\lf\cha_s,\l\rf)=\exp(2\pi\sqrt{-1}\lf\cha_s,\l\rf)$
or that $(-1+(q+1)/2)\lf\cha_s,\l\rf=0$, or that $2t\lf\cha_s,\l\rf=0$. This follows from 0.5(b).

\mpb

We show:

(c) If $sw=ws$, $|sw|<|w|$ then $\d_{w,\l;s}$ defined as in 0.3 in terms of $q$ is equal to $1$.
\nl
It is enough to show that 
$$\exp(2\pi\sqrt{-1}((q-1)/2)\lf\cha_s,\l\rf)=1$$ or that $((q-1)/2)\lf\cha_s,\l\rf=0$.
 Since $\l\in \bX_{q'}$ we have $(q'-1)\lf\cha_s,\l\rf=0$ by the argument at the end of 0.3.
We have $(q-1)/2=(q'-1)(q'+1)/2$ where $q'+1\in2\ZZ$ hence
$$((q-1)/2)\lf\cha_s,\l\rf=((q'+1)/2)(q'-1)\lf\cha_s,\l\rf=0.$$
 This proves (c).

\proclaim{Proposition 4.3} Let $\qq$ be an indeterminate and let $\ti\MM$ denotes the
$\CC(\qq)$-vector space with basis $\{\ti\aa_{w,\l};(w,\l)\in\tX_1\}$.
There is a unique action of the braid group of $W$ on $\ti\MM$
in which the generators $\{\ct_s;s\in S\}$ of the braid group applied to the basis elements
of $\ti\MM$ are as follows. (We write $\D=1$ if $s\in W_\l$ and $\D=0$ if $s\n W_\l$.)

(a) $\ct_s\ti\aa_{w,\l}=\ti\aa_{sws,\l}$ if $sw\ne ws,|sw|>|w|$;

(b) $\ct_s\ti\aa_{w,\l}=\ti\aa_{sws,s\l}+\D(\qq-\qq\i)\ti\aa_{w,\l}$ if $sw\ne ws,|sw|<|w|$;

(c) $\ct_s\ti\aa_{w,\l}=\d'_{w,s\l;s}\ti\aa_{w,s\l}+\D(\qq+1)\ti\aa_{sw,\l}$ if $sw=ws,|sw|>|w|$;

(d) $\ct_s\ti\aa_{w,\l}=\D(1-\qq\i)\ti\aa_{sw,\l}+\D(\qq-\qq\i)\ti\aa_{w,\l}-\ti\aa_{w,s\l}$ if $sw=ws,|sw|<|w|$.
\nl
Here $\d'_{w,s\l;s}=\pm1$ is as in 0.5. (It is $1$ in the simply laced case; it is also $1$ if $\D=1$.)
\endproclaim
It is enough to prove the proposition with $\ti\MM$ replaced by 
the $\CC(\qq)$-vector space $\ti\MM_\co$ with basis $\{\ti\aa_{w,\l};(w,\l)\in\tX_1,\l\in\co\}$ 
where $\co$ is any $W$-orbit in $\bX$.

Assume first that $\co$ is as in 4.1 and let $\fe,\fQ,\tfQ$ be as in 4.2. Let $\tfQ'=\{q\in\tfQ; 2\fe<q^2-1\}$. 
Clearly, $\tfQ'$ is an infinite set.

Let $M_{\co}$ be the $\CC$-vector space with basis 
$\{\ti\aa_{w,\l};(w,\l)\in \tX_1;\l\in\co\}$.
By 4.1(b) we can identify $M_\co$ with a subspace of $M_q$ (for any $q\in\fQ$) by
$\ti\aa_{w,\l}\m a_{w,\l}$. This subspace of $M_q$ is stable under the operators
$\ct_s,s\in S$ attached in 0.4 to $z=\fe$, provided that $q\in\tfQ'$. (Note for $q\in\tfQ'$
we have $2z\n (q^2-1)\ZZ$ since $0<2fe<q^2-1$.) Hence $\ct_s:M_q@>>>M_q$ can be regarded as an operator 
$\ct_s^{(q)}:M_\co@>>>M_\co$ for any $q\in\tfQ'$. This operator is given by matrix in the basis of $M_\co$ is given by
Laurent polynomials in $q$ with integer coefficients independent of $q$. (This follows from the formulas 0.4(a)-(h), from
4.2(a),(b),(c) and from the equality $\exp(2\pi\sqrt{-1}(w:s)\fe\la\cha_s,\l\ra)=1$ for $\l\in\co$.) 
 Since $q$ runs through an infinite set, we 
deduce that the braid group relations satisfied by the $\ct_s^{(q)}$ remain valid
when $q$ is replaced by the indeterminate $\qq$. We see that if we identify $\ti\MM_\co=\CC(\qq)\ot M_\co$, then
there is a unique action of the braid group of $W$ on $\ti\MM_\co$
in which the generators $\{\ct_s;s\in S\}$ of the braid group applied to the basis elements
of $\ti\MM_\co$ are as in (a)-(d) above. 

We now consider a $W$-orbit $\co$ in $\bX$ which is not necessarily as in 4.1.
We choose $\x_0\in X_K$ such that the image of $x_0$ in $\bX$ belongs to $\co$.
Let $\fH$ be the collection of affine hyperplanes 

$\{\x\in X_K; \la\cha,\x\ra=e\}$ for various $\cha\in\che R,e\in\ZZ$;

$\{\x\in X_K; w(\x)=\x+x\}$ for various $w\in W-\{1\}, x\in X$;

$\{\x\in X_K; w(\x)=-\x+x\}$ for various $w\in W_2, x\in X$ such that $w+1$ is not identically zero on $X$.
\nl
We can find $\x'_0\in X_\QQ$ such that a hyperplane in $\fH$ contains $\x_0$ if and only if it contains
$\x'_0$. Let $\co'$ be the $W$-orbit of the image of $\x'_0$ in $\bX$.
There is a unique $W$-equivariant bijection $j:\co'@>\si>>\co$ under which the image of $\x'_0$ in $\bX$ corresponds to 
the image of $\x_0$ in $\bX$.
We define an isomorphism $\ti\MM_{\co'}@>\si>>\ti\MM_\co$ by $\ti\aa_{w,\l'}\m\ti\aa_{w,j(\l')}$. This isomorphism is
compatible with the operators $\ct_s$ on these two vector spaces. Since there operators satisfy the braid group
relations on $\ti\MM_\co'$ (by the first part of the proof) they will satisfy the braid group relations on $\ti\MM_\co$.
This completes the proof of the proposition.

\subhead 4.4\endsubhead
Let $v$ be an indeterminate such that $v^2=\qq$. Let
$\MM=\CC(v)\ot_{\CC(\qq)}\ti\MM$.
We consider the basis $\{\aa_{w,\l};(w,\l)\in\tX_1\}$ defined by
$\aa_{w,\l}=v^{||w||}\ti\aa_{w,\l}$
where $||w||$ is as in 2.12. The linear maps $\ct_s:\ti\MM@>>>\ti\MM$ with $s\in S$
 extend to linear maps $\ct_s:\MM@>>>\MM$ which satisfy the equalities in
 0.6. Thus Theorem 0.6 is a consequence of Proposition 4.3.

\subhead 4.5\endsubhead
Let $\HH$ be the $\CC(v)$-vector space with basis $\{\ct_{w,\l};(w,\l)\in W\bX\}$.
There is a unique structure of associative $\CC(v)$-algebra (without $1$ in general) on $\HH$
such that (a),(b) below hold.
$$\ct_{w,\l}\ct_{w',\l'}=\d_{w\i(\l),\l'}\ct_{ww',\l'}\tag a$$ 
if $(w,\l)\in W\bX,(w',\l')\in W\bX$, $|ww'|=|w|+|w'|$;
$$\ct_{s,\l}\ct_{s,\l'}=\d_{\l,\l'}\ct_{1,\l'}+\D\d_{s(\l),\l'}(v^2-v^{-2})\ct_{s,\l'}\tag b$$
if $s\in S,\l\in\bX,\l'\in\bX$ (here $\D=1$ if $s\in W_\l$ and $\D=0$ if $s\n W_\l$).
We call $\HH$ the {\it extended Hecke algebra}. 
This algebra has been studied in \cite{\CDGVII}, \cite{\MND} (at least when $K=\QQ$). It is similar but not the same to an
algebra studied in \cite{\MS}.

For any $w\in W$ we define a linear map $\ct_w:\MM@>>>\MM$ by 
$\ct_w=\ct_{s_1}\ct_{s_2}\do\ct_{s_k}$
where $s_1,s_2,\do,s_k$ are elements of $S$ such that $w=s_1s_2\do s_k$, $|w|=k$. By 0.6, 
this is independent
of the choice of $s_1,\do,s_k$. For $\l\in\bX$ we define a linear map $1_\l:\MM@>>>\MM$ by
$1_\l(\aa_{w,\l'})=\d_{\l,\l'}\aa_{w,\l'}$ for any $(w,\l')\in\tX_1$. 
For $(w,\l)\in W\bX$ we define a linear map $\ct_{w,\l}:\MM@>>>\MM$ as the composition $\ct_w1_\l$.
These maps define an $\HH$-module structure on $\MM$. (This follows from 0.6; the relation
(b) on $\MM$ can be deduced from the analogous relation in $M_q$.)
From (b) we deduce that $\ct_s\i:\MM@>>>\MM$ is well defined and we have
$$\ct_s\i=\ct_s-(v^2-v^2)\i\sum_{\l\in\bX;s\in W_\l}1_\l.\tag c$$
(The last sum may be infinite but at most  one term in the sum applied to a given basis 
element of $\MM$ can be nonzero.) It follows that for any $w\in W$, $\ct_w:\MM@>>>\MM$ 
is invertible. Its inverse satisfies $\ct_{w_1w_2}\i=\ct_{w_2}\i\ct_{w_1}\i:\MM@>>>\MM$
for any $w_1,w_2$ in $W$ such that $|w_1w_2|=|w_1+|w_2|$.

For any $W$-orbit $\co$ in $\bX$ we denote by $\HH_\co$ the subspace of $\HH$ spanned by 
$$\{\ct_{w,\l};(w,\l)\in W\T\co\}.$$
This is a subalgebra of $\HH$, this time with unit, namely $\sum_{\l\in\co}\ct_{1,\l}$. 

For any $w\in W$ we set $\ct_w=\sum_{\l\in\co}\ct_{w,\l}\in\HH_\co$; for any $\l\in\co$ we set
$1_\l=\ct_{1,\l}\in\HH_\co$. We see that the elements $\ct_w,1_\l$ exist separately in $\HH_\co$, not only in
the combination $\ct_{w,\l}=\ct_w1_\l$.

We denote by $\MM_\co$ the subspace of $\MM$ spanned by 
$\{\aa_{w,\l};(w,\l)\in\tX_1,\l\in\co\}$. Note that the $\HH$-module structure on $\MM$
restricts to an $\HH_\co$-module structure on $\MM_\co$.

\head 5. On the structure of the $\HH$-module $\MM$\endhead
\subhead 5.1\endsubhead
In this section we assume that $\kk=\CC$. 
For $\l\in\bX$ let $\che R_\l=\{\cha\in\che R;\lf\cha,\l\rf=0$, 
$\che R_\l^+=\che R_\l\cap\che R^+$. Then $\cha R_\l$ is the set of coroots of a root system and $\che R_\l^+$
is a set of positive coroots for it. Let $\che R_\l^-=\che R_\l-\che R_\l^+$.
Let $\che\Pi_\l$ be the set of simple coroots for $\che R_\l$ contained in
$\che R^+_\l$. For each $\b\in\che R$ let $s_\b:Y@>>>Y$ be the reflection in $W$ such that $s_\b(\b)=-\b$.
Let $W_\l$ be the subgroup of $W$ generated by $\{s_\b;\b\in\che R_\l\}$. This is a Coxeter group with generators 
$\{s_\b;\b\in\che\Pi_\l\}$ and with length function $w\m|w|_\l=\sha(\b\in\che R_\l^+;w(\b)\in\che R_\l^-)$.
Note that for $s\in S$ the condition that $s\in W_\l$ coincides with the condition denoted in the
same way in 0.1; this follows from \cite{\MND, 1.2(c)}.

If $w\in W$ then there is a unique element $z\in wW_\l$ such that $z(\che R_\l^+)\sub\che R^+$; we have
$|z|<|zu|$ for any $u\in W_\l-\{1\}$; we write $z=\min(wW_\l)$. (See \cite{\MND, 1.2(e)}).

We now fix an integer $m\ge1$. We fix a $W$-orbit $\co$ in $\bX_m$. For any $\l,\l'$ in $\co$ we set 
$$[\l',\l]=\{z\in W;\l'=z(\l),z=\min(zW_\l)\}=\{ z\in W;\l'=z(\l),z(\che R_\l^+)=\che R_{\l'}^+\}.$$
Clearly, 

(a) $[\l,\l']=[\l',\l]\i$; moreover, if $\l,\l',\l''$ are in $\co$ then $[\l'',\l'][\l',\l]\sub[\l'',\l]$.
\nl
Hence the group structure on $W$ makes 

(b) $\Xi:=\{(\l',z,\l)\in\co\T W\T\co; z\in[\l',\l]\}$
\nl
into a groupoid, see \cite{\MND, 1.2(f)}.

\subhead 5.2\endsubhead
If $\l\in\bX$ then $\cha R_\l\sub\cha R_{-m\l}$. If $(w,\l)\in\tX_m$ then 
$\sha(\che R_\l)=\sha(\che R_{-m\l})$ so that $\che R_\l=\che R_{-m\l}$ and $W_\l=W_{-m\l}$. We show:

(a) {\it If $\l\in\bX_m$ and $z\in[-m\l, m]$ then $z(\che R_\l^+)=\che R_\l^+$ so that
$\io_z:u\m zuz\i$ is a Coxeter group automorphism of $W_\l$.}
\nl
We have $z(\che R_\l)=\che R_{z\l}=\che R_{-m\l}=\che R_\l$; 
moreover since $z(\che R_\l^+)\sub \che R^+$ we have $z(\che R_\l^+)=\che R_\l^+$.
This proves (a).

\mpb

Let $\tX_m^0=\{(z,\l)\in W_2\T\bX; z\in[-m\l,\l]\}$. Note that $\tX_m^0\sub\tX_m$.
For $(z,\l)\in\tX_m^0$ let $I_{(z,\l)}=\{u\in W_\l; \io_z(u)u=1\}$ be the 
set of $\io_z$-twisted involutions of $W_\l$.
If $u\in I_{z,\l}$ then $(zu,\l)\in\tX_m$; indeed we have $(zu)^2=1$ and $zu(\l)=z(\l)=-m\l$. Conversely, 

(b) if $(w,\l)\in\tX_m$ we have $(w,\l)=(zu,\l)$
for a well defined $(z,\l)\in\tX_m^0$ and $u\in I_{z,\l}$.
\nl
 Indeed, let $z=\min(wW_\l)$. Since $w(l)=-m\l$ we have also $z(\l)=-m\l$ hence $z\in[-m\l,\l]$.
We have $w=zu$ where $u\in W_\l$. We have $w=w\i=u\i z\i=z\i zuz\i=z\i\io_z(u)$.
Since $\io_z(u)\in W_\l$ (see (a)) we have $w\in z\i W_\l$.
Since $z(\che R_\l^+)=\che R_\l^+$ we must have also $z\i(\che R_\l^+)=\che R_\l^+$
so that $z\i=\min(wW_\l)$. It follows that $z=z\i$ so that $(z,\l)\in\tX_m^0$.
Since $1=w^2=(zu)^2$ we see that $\io_z(u)u=1$ so that $u\in I_{z,\l}$. This proves (b).

We see that

(c) we have a bijection $\sqc_{(z,\l)\in\tX_m^0}I_{z,\l}@>\si>>\tX_m$ given
by $(z,\l,u)\m (zu,\l)$ where $(z,\l)\in\tX_m^0,u\in I_{z,\l}$.

\subhead 5.3\endsubhead
Let $\Xi$ be as in  5.1(b). Let $\Xi^0=\{(z,\l)\in\tX_m^0;\l\in\co\}$.

We can view  $\Xi^0_m$ as a subset of $\Xi$ by $(z,\l)\m(-m\l,z,\l)$.
This subset is the fixed point set of the antiautomorphism
$$(\l',z,\l)\m(\l',z,\l)^*:=(-m\l,z\i,-m\l')$$
 of the groupoid $\Xi$ (the composition of the inversion 
$(\l',z,\l)\m(\l,z\i,\l')$ with the involutive automorphism
$(\l',z,\l)\m(-m\l',z,-m\l)$ of the groupoid $\Xi$). Hence this subset
can be viewed as the set of $*$-twisted "involutions" of this groupoid.

\mpb

Until the end of 5.8 we assume that $m=1$. From 0.6 we deduce

(a) If $(w,\l)\in\tX_1$, $s\in S$ and $s\n W_\l$ then $\ct_s(\aa_{w,\l})=\pm\aa_{sws,s\l}$.
\nl
Note also that in $\HH_\co$, for $s\in S,w\in W,\l\in\co$ we have

(b) $\ct_s\ct_w1_\l=\ct_{sw}1_\l$ if $s\n W_{w(\l)}$; $\ct_w\ct_s1_\l=\ct_{ws}1_\l$ if $s\n W_\l$.

\proclaim{Lemma 5.4}Let $\l\in\co$. Let $(w,\l)\in\tX_1$, $z\in[\l,\l]$. Then
$(zwz\i,\l)\in\tX_1$ and $\ct_z\aa_{w,\l}=\pm \aa_{zwz\i,\l}$.
\endproclaim
The proof is similar to that of \cite{\MND, 1.4(c)}.
We have $w(\l)=-\l$ hence $zwz\i(\l)=-\l$ since $z(\l)=\l$. Thus
$(zwz\i,\l)\in\tX_1$.

We write $z=s_ks_{k-1}\do s_1$ where $s_1,\do,s_k$ are in $S$, $|z|=k$.
As in the proof of \cite{\MND, 1.4(c)} we have $s_1\n W_\l$,
$s_1s_2s_1\n W_\l$, $\do$ $s_1s_2\do s_k\do s_2s_1\n W_\l$.
We have $\ct_{s_1}\aa_{w,\l}=\pm \aa_{s_1ws_1},s_1\l)$ since $s_1\n W_\l$, see
5.3(a). We have  $\ct_{s_2}\aa_{s_1ws_1,s_1\l}=\pm\aa_{s_2s_1ws_1s_2,s_2s_1\l}$
since $s_2\n W_{s_1\l}$, see 5.3(a). Continuing in this way we get
$$\ct_{s_k}\aa_{s_{k-1}\do s_1ws_1\do s_{k-1},s_{k-1}\do s_1\l}
 =\pm\aa_{s_k\do s_1ws_1\do s_k,s_k\do s_1\l}.$$
Combining these equalities we get
$$\ct_z\aa_{w,\l}=\ct_{s_k}\do\ct_{s_1}\aa_{w,\l}=
\pm \aa_{s_k\do s_1ws_1\do s_k,s_k\do s_1\l}=\pm\aa_{zwz\i,z\l}
=\pm\aa_{zwz\i,\l}.$$
The lemma is proved.

The following result is a generalization of the lemma above.

\proclaim{Lemma 5.5}Let $(w,\l)\in\tX_1$, $z\in[\l',\l]$ where $\l,\l'$ are in
$\co$. Then $(zwz\i,\l')\in\tX_1$ and $\ct_z\aa_{w,\l}=\pm \aa_{zwz\i,\l'}$.
\endproclaim
The proof is similar to that of \cite{\MND, 1.4(d)}.
We have $w(\l)=-\l$ hence $zwz\i(\l')=-\l'$ since $z\i(\l')=\l'$.
Thus $(zwz\i,\l')\in\tX_1$.

Since $\l,\l'$ are in the same $W$-orbit, we can find $r\ge0$ and
$s_1,s_2,\do,s_r$ in $S$ such that, setting 
$$\l_0=\l,\l_1=s_1\l,\l_2=s_2s_1\l,\do,\l_r=s_r\do s_2s_1\l,$$
we have $\l_0\ne\l_1\ne\l_2\ne\do\ne\l_r=\l'$.
For $j=1,\do,r$, we have $s_j\n W_{\l_{j-1}}$  since
$s_j(\l_{j-1}=\l_j\ne\l_{j-1}$' hence $s_j$ has minimal length in
$s_jW_{\l_{j-1}}$ and $s_j\in[\l_j,\l_{j-1}]$.
It follows that $s_r\do s_2s_1\in[\l_r,\l_0]=[\l',\l]$ (we use 5.1(a)).
We define $\tz\in W$ by $z=s_r\do s_2s_1\tz$. Then $\tz\in[\l,\l]$ (we use
again 5.1(a)).
For $j\in[1,r]$ we have $s_j\n W_{s_{j-1}\do s_1\l}$ (since $\l_j\ne\l_{j-1}$)
hence, using 5.3(a) we have
$$\ct_{s_j}\aa_{s_{j-1}\do s_1\tz w\tz\i s_1\do s_{j-1},s_{j-1}\do s_1\l}=
\pm\aa_{s_js_{j-1}\do s_1\tz w\tz\i s_1\do s_{j-1}s_j,s_js_{j-1}\do s_1\l}.$$
Applying this repeatedly we deduce
$$\ct_{s_r}\do\ct_{s_2}\ct_{s_1}\aa_{\tz w\tz\i,\tz\l}
=\pm\aa_{s_r\do s_2s_1\tz w\tz\i s_1s_2\do s_r,s_r\do s_2s_1\tz\l}
=\pm\aa_{zwz\i,z\l}.$$
We now apply 5.4 with $z$ replaced by $\tz$; we see that
$\ct_{\tz}\aa_{w,\l}=\pm \aa_{\tz w\tz\i,\l}$. Substituting this in the
previous equation we obtain
$$\ct_{s_r}\do\ct_{s_2}\ct_{s_1}\ct_{\tz}\aa_{w,\l}=\pm\aa_{zwz\i,z\l}.\tag a$$
For $j\in[1,r]$ we have $s_j\n W_{s_{j-1}\do s_1\l}$ (as above)
hence, using 5.3(b) we have
$$\ct_{s_j}\ct_{s_{j-1}\do s_1\tz}\aa_{w,\l}
=\ct_{s_js_{j-1}\do s_1\tz}\aa_{w,\l}.$$
Applying this repeatedly we deduce
$$\ct_{s_r}\do\ct_{s_2}\ct_{s_1}\ct_{\tz}\aa_{w,\l}=
\ct_{s_r\do s_2s_1\tz}\aa_{w,\l}=\ct_z\aa_{w,\l}.$$
Combining this with (a) gives
$$\ct_z\aa_{w,\l}=\pm\aa_{zwz\i,z\l}.$$
The lemma is proved.

\proclaim{Lemma 5.6} Let $(z,\l)\in\Xi^0$ and let $u\in W_\l$. Let
$\a\in\che\Pi_\l$. We set $\s=\s_\a$; note that $|\s|_\l=1$. Recall that $u\m\io_z(u)=zuz\i$ is an 
involutive Coxeter group automorphism of $W_\l$. For any $u\in W_\l$ we have

(a) $\ct_\s\aa_{zu,\l}=e_1\aa_{z\io_z(\s)u\s,\l}$ if $u\s\ne\io_z(\s)u,|u\s|_\l>|u|_\l$;

(b) $\ct_\s\aa_{zu,\l}=e_2\aa_{z\io_z(s)u\s,\l}+e_3(v^2-v^{-2})\aa_{zu,\l}$ if $u\s\ne\io_z(\s)u,|u\s|_\l<|u|_\l$;

(c) $\ct_\s\aa_{zu,\l}=e_4\aa_{zu,\l}+e_5(v+v\i)\aa_{zu\s,\l}$ if $u\s=\io_z(\s)u$, $|u\s|_\l>|u|_\l$;

(d)  $\ct_\s\aa_{zu,\l}=e_6(v-v\i)\aa_{zu\s,\l}+e_7(v^2-v^{-2}-1)\aa_{zu,\l}$ if $u\s=\io_z(\s)u$, $|u\s|_\l<|u|_\l$,
\nl
where $e_1,\do,e_7\in\{1,-1\}$.
\endproclaim
As in the proof of \cite{\MND,1.4(f)} we can find $s_1,s_2,\do,s_r$ in $S$ such that
$\s=s_1s_2\do s_{r-1}s_rs_{r-1}\do s_1$, $|\s_\a|=2r-1$,
$s_1s_2\do s_{j-1}s_js_{j-1}\do s_1\n W_\l$ for $j=1,2,\do r-1$.
We argue by induction on $r\ge1$. When $r=1$ the result follows from 0.6.
(Note that $z\io_z(\s)u\s=\s zu\s$, the condition $u\s=\io_z(\s)u$ is equivalent to $zu\s=\s zu$ and if
$|\s|=1$ the condition $|u\s|_\l>|u|_\l$ is equivalent to $|u\s|>|u|$.)
Assume now that $r\ge2$. We set $s=s_1$, $\l'=s\l$, $\b=s(\a)\in R^+_{\l'}$, $u'=sus$, $z'=szs$, $\s'=s_\b=s\s s$.
We have $(z',\l')\in\Xi_\co^0$, $u'\in W_{\l'}$ and $\s'\in W_{\l'}$,
$|\s'|_{\l'}=1$, $|\s'|=|\s|-2$. Moreover, we have $s\n W_\l$. By the induction hypothesis we have

(a${}'$) $\ct_{\s'}\aa_{z'u',\l'}=e'_1\aa_{\s'z'u'\s',\l'}$ if $u'\s'\ne z'\s'z'u',|u'\s'|_{\l'}>|u'|_{\l'}$;

(b${}'$)  $\ct_{\s'}\aa_{z'u',\l'}=e'_2\aa_{\s'z'u'\s',\l'}+e'_3(v^2-v^{-2})\aa_{z'u',\l'}$ if $u'\s'\ne z'\s'z'u',
|u'\s'|_{\l'}<|u'|_{\l'}$;

(c${}'$) $\ct_{\s'}\aa_{z'u',\l'}=e'_4\aa_{z'u',\l'}+e'_5(v+v\i)\aa_{z'u'\s',\l'}$ if $u'\s'=z'\s'z'u',
|u'\s'|_{\l'}>|u'|_{\l'}$;

(d${}'$) $\ct_{\s'}\aa_{z'u',\l'}=e'_6(v-v\i)\aa_{z'u'\s',\l'}+e'_7(v^2-v^{-2}-1)\aa_{z'u',\l'}$ if 
$u'\s'=z'\s'z'u', |u'\s'|_{\l'}<|u'|_{\l'}$,
\nl
where $e'_1,\do,e'_7\in\{1,-1\}$. By 5.3(a), 5.3(b) we have 
$$\ct_s\ct_{\s'}\aa_{z'u',\l'}=\ct_\s\ct_s\aa_{z'u',\l'}=\pm\ct_\s\aa_{zu,\l}.$$
Moreover, by 5.3(a) we have
$$\ct_s\aa_{z'u',\l'}=\aa_{zu,\l}, \ct_s\aa_{\s'z'u'\s',\l'}=\aa_{\s zu\s,\l},$$
$\ct_s\aa_{z'u'\s',\l'}=\aa_{zu\s,\l}$.
Hence (a)-(d) for $\s,z,u$ follow from (a${}'$)-(d${}'$) by applying $\ct_s$ to both sides.
Here we use that the condition that $z'u'\s'=\s'z'u'$ is equivalent to the condition $zu\s=\s zu$ 
and the inequality $|u'\s'|_{\l'}>|u'|_{\l'}$ is equivalent to the inequality
$|u\s|_\l>|u|_\l$ (conjugation by $s$ is a Coxeter group isomorphism $W_{\l'}@>>>W_\l$).
The lemma is proved.

\subhead 5.7\endsubhead
For any $\l\in\co$ let $\HH_\l$ be the $\CC(v)$-subspace of $\HH_\co$ spanned by $\{\ct_u1_\l;u\in W_\l\}$.
This is a subalgebra of $\HH_\co$ with unit $1_\l$; it can be identified with the
Hecke algebra of the Coxeter group $W_\l$ (see \cite{\MND, 1.4(g),(h)})
so that the standard generators of the
last algebra correspond to the elements $\ct_{s_\a}1_\l$ of $\HH_\l$ with $\a\in\che\Pi_\l$.

For $(z,\l)\in\Xi^0$ let $\MM_{z,\l}$ be the subspace of $\MM$ spanned by $\{\aa_{zu,\l};u\in I_{z,\l}\}$.
From 5.6 we see that $\MM_{z,\l}$ is an $\HH_\l$-module and that the action of the generators of $\HH_\l$ on 
$\MM_{z,\l}$ is given by a formula which is the same (except for the apparition of certain signs $e_j$)
as the formula for the action of the generators of the Hecke
algebra of $W_\l$ on the module based on the twisted involutions in $W_\l$ constructed in \cite{\LV}.

\subhead 5.8\endsubhead
We have a direct sum decomposition $\HH_\co=\op_{(\l',z,\l)\in\Xi}\ct_z\HH_\l$; moreover,
$\{\ct_z1_\l\ct_u1_\l;(\l',z,\l)\in\Xi, u\in W_\l\}$ is a basis of $\HH_\co$ compatible with this
decomposition and it coincides with the the basis $\{\ct_w1_\l;(w,\l)\in\tX_1,\l\in\co\}$ of $\HH_\co$.
(See \cite{\MND, 1.4(d)}.) 
Similarly, by 5.2(b), we have a direct sum decomposition $\MM_\co=\op_{(\tz,\ti\l)\in\Xi^0}\MM_{\tz,\ti\l}$
where $\MM_{\ti z,\ti\l}$ is as in 5.7.
From Lemmas 5.5, 5.6 we see that the direct sums decompositions
of $\HH_\co$ and $\MM_\co$ are compatible in the following sense:
$$(\ct_z\HH_\l)\MM_{\ti z,\ti\l}\sub \d_{\ti\l,\l}\MM_{z\ti zz\i, z(\ti\l)}.$$
Moreover the action of the basis element $\ct_z\ct_u1_\l=(\ct_z1_\l)(\ct_u1_\l)$ of $\HH_\co$  on a basis element
$\aa_{\tz u',\ti\l}$ of $\MM_\co$ is particularly simple:
the operator $\ct_z1_\l$ applied to a basis element $\aa_{\tz u',\ti\l}$ is $\pm\d_{\ti\l,\l}$  times another basis element;
the operator $\ct_u1_\l$ applied to a basis element  $\aa_{\tz u',\ti\l}$ is as in 5.7 if $\ti\l=\l$ and is zero if $\ti\l\ne\l$..

\subhead 5.9\endsubhead
Results similar to those in 5.4-5.8 hold for $M_\co$ when $m=q$ with $(p,q)$ as in 0.2 and $\co\sub \bX_m$ except that in this case
the $\pm$ signs in 5.4-5.8 have to be replaced by roots of $1$ of possibly higher order.

\head 6. Proof of Theorem 0.9\endhead
\subhead 6.1\endsubhead
We now fix an integer $m\ge1$. Recall from 0.8 that $\MM_m$ is the $\CC)v)$-vector space with basis
 $\{\aa_{w,\l};(w,\l)\in\tX_m\}$.
We fix a $W$-orbit $\co$ in $\bX_m$. Let $\MM_\co$ be the subspace of $\MM_m$ spanned by
$\{\aa_{w,\l};(w,\l)\in\tX_m,\l\in\co\}$.

For any $\l\in\co$ let $\HH_\l$ be as in 5.7.
For $(z,\l)\in\Xi^0$ let $\MM_{z,\l}$ be the subspace of $\MM_\co$ spanned by $\{\aa_{zu,\l};u\in I_{z,\l}\}$.
By \cite{\LV} applied to the Coxeter group $W_\l$ with the involutive automorphism $\io_z$, there is a well defined 
$\HH_\l$-module structure $(h,\x)\m h\cir\x$ on $\MM_{z,\l}$ such that
for any $u\in W_\l$ and any $\s=s_\a$, $\a\in\che\Pi_\l$ we have

(a) $(\ct_\s1_\l)\cir\aa_{zu,\l}=\aa_{z\io_z(\s)u\s,\l}$ if $u\s\ne\io_z(\s)u,|u\s|_\l>|u|_\l$;

(b) $(\ct_\s1_\l)\cir\aa_{zu,\l}=\aa_{z\io_z(s)u\s,\l}+(v^2-v^{-2})\aa_{zu,\l}$ if $u\s\ne\io_z(\s)u,|u\s|_\l<|u|_\l$;

(c) $(\ct_\s1_\l)\cir\aa_{zu,\l}=\aa_{zu,\l}+(v+v\i)\aa_{zu\s,\l}$ if $u\s=\io_z(\s)u$, $|u\s|_\l>|u|_\l$;

(d)  $(\ct_\s1_\l)\cir\aa_{zu,\l}=(v-v\i)\aa_{zu\s,\l}+(v^2-v^{-2}-1)\aa_{zu,\l}$ if $u\s=\io_z(\s)u$, $|u\s|_\l<|u|_\l$.

\subhead 6.2\endsubhead
By \cite{\MND, 1.4(d)} the basis $\{\ct_w1_\l;(w,\l)\in\tX_1,\l\in\co\}$ of $\HH_\co$ coincides with
$\{\ct_u\ct_z1_\l;(\l',z,\l)\in\Xi, u\in W_{\l'}\}$. We define a bilinear multiplication $\HH_\co\T\MM_\co@>>>\MM_\co$ (denoted by
$(h,\x)\m h\bul\x$) by the rule
$$(\ct_u\ct_z1_\l)\bul\aa_{\tz\ti u,\ti\l}=0$$
if $\l\ne\ti\l$, while if $\l=\ti\l$,
$$(\ct_u\ct_z1_\l)\bul\aa_{\tz\tu,\ti\l}=(\ct_u1_{\l'})\cir\aa_{(z\tz z\i)(z\tu z\i),\l'}$$
for $(\l',z,\l)\in\Xi,u\in W_{\l'},(\tz,\ti\l)\in\Xi^0,\tu\in W_{\ti\l}$, where $\cir$ is as in 6.1 with $\l$ replaced by 
$\l'$. (We have $(z\tz z\i,\l')\in\Xi^0$ and $z\tu z\i\in W_{\l'}$.) We show:

(a) {\it this is an $\HH_\co$-module structure.}
\nl
It is enough to show that for 
$$(\l',z,\l)\in\Xi,(\l'_1,z_1,\l_1)\in\Xi, u\in W_{\l'},u_1\in W_{\l'_1},(\tz,\ti\l)\in\Xi^0,\tu\in W_{\ti\l},$$
 with $\l'=\l_1,\l=\ti\l$ we have
$$(\ct_{u_1}\ct_{z_1}1_{\l_1})\bul(\ct_u\ct_z1_\l)\bul\aa_{\tz\tu,\ti\l}=
(\ct_{u_1}\ct_{z_1uz_1\i}\ct_{z_1z}1_\l)\bul\aa_{\tz\tu,\ti\l}$$
or that
$$\align&(\ct_{u_1}1_{\l'_1})\cir((\ct_{z_1uz_1\i}\ct_{z_1zz_1\i}1_{z_1\l})\bul\aa_{(z_1\tz z_1\i)(z_1\ti uz_1\i),z_1\l})\\&=
\sum_{u_2\in W_{\l'_1}}(\ct_{u_2}\ct_{z_1z}1_\l)\bul\aa_{\tz\ti u,\ti\l}\endalign$$
where we have written $\ct_{u_1}\ct_{z_1uz_1\i}1_{\l'_1}=\sum_{u_2\in W_{\l'_1}}\g_{u_2}\ct_{u_2}1_{\l'_1}$, $\g_{u_2}\in\CC(v)$.
(We have used \cite{\MND, 1.4(d),(e)}.) We have
$$\align&(\ct_{z_1uz_1\i}\ct_{z_1zz_1\i}1_{z_1\l})\bul\aa_{z_1\tz z_1\i)(z_1\ti uz_1\i),z_1\l})\\&=
(\ct_{z_1uz_1\i}1_{\l'_1})\cir\aa_{(z_1z\tz z\i z_1\i)(z_1z\ti uz\i z_1\i),\l'_1}).\endalign$$
We have
$$\align&\sum_{u_2\in W_{\l'_1}}(\ct_{u_2}\ct_{z_1z}1_\l)\bul\aa_{\tz\ti u,\ti\l}\\&=
\sum_{u_2\in W_{\l'_1}}(\ct_{u_2}1_{\l'_1})\cir \aa_{(z_1z\tz z\i z_1\i)(z_1z\tu z\i z_1\i),\l'_1}.\endalign$$
Thus it is enough to prove
$$\align&(\ct_{u_1}1_{\l'_1})\cir((\ct_{z_1uz_1\i}1_{\l'_1})\cir\aa_{(z_1z\tz z\i z_1\i)(z_1z\tu z\i z_1\i),\l'_1})\\&
=\sum_{u_2\in W_{\l'_1}}(\ct_{u_2}1_{\l'_1})        \cir \aa_{(z_1z\tz z\i z_1\i)(z_1z\tu z\i z_1\i),\l'_1}.\endalign$$
This follows from the fact that $\cir$ defines a module structure. This proves (a).

\subhead 6.3\endsubhead
Let $\HH_m$ be the $\CC(v)$-vector space with basis $\{\ct_{w,\l};(w,\l)\in W\T\bX_m\}$.
Note that $\HH_m$ is a subalgebra of $\HH$.
There is a unique $\HH_m$-module structure $(h,\x)\m h\bul\x$ on $\MM_m$ (see 0.8) such that 
for any two orbits $\co,\co'$ in $\bX_m$ and any $h\in\HH_\co,\x\in\MM_{\co'}$ we have
$h\bul\x=0$ if $\co\ne\co'$ and $h\bul\x$ is as in 6.2 if $\co=\co'$.

\subhead 6.4\endsubhead
We now prove Theorem 0.9.
It is enough to show that 0.9(a)-(b) hold when $\ct_s$ is replaced by $\ct_s1_\l\in\HH_m$ acting on $\MM_m$
via the $\HH_m$-module structure on $\MM_m$. We can write $w=zu$ where $(z,\l)\in\Xi^0$ and $u\in W_\l$.
If $s\in W_\l$ then $s=\s$ as in 6.1 and the desired formulas
follow from 6.1. If $s\n W_\l$ then $s$ has minimal length in $sW_\l$ hence $s\in[s(\l),\l]$.
Then by definition we have $(\ct_s1_\l)\bul\aa_{w,\l}=\aa_{sws,s\l}$ and the desired formulas hold again.
This proves Theorem 0.9.

\subhead 6.5\endsubhead
In \cite{\MND}, an affine analogue of $\HH$ is considered; it has a basis indexed by the semidirect product
$\tW\bX$ where $\tW$ is an affine Weyl group acting on $\bX$ via its quotient $W$. The analogue of Theorem 0.9
continues to hold in this case (with the same proof).

\head 7. Bar operator\endhead
\subhead 7.1\endsubhead
Let $m$ be an integer $\ge1$.
In this section we construct a bar operator on $\MM_m$ generalizing a definition in \cite{\LV}. 
To do this we will use the method of \cite{\BAR}.

For $s\in S$ the operator $\ct_s:\MM_m@>>>\MM_m$ in 0.9 has an inverse $\ct_s\i$. For $w\in W$ we set
$\ct_w=\ct_{s_1}\do\ct_{s_k}:\MM_m@>>>\MM_m$, $\ct_w\i=\ct_{s_k}\i\do\ct_{s_1}\i:\MM_m@>>>\MM_m$, 
where $w=s_1s_2\do s_k$ with $s_1,\do,s_k$ in $S$, $|w|=k$.

Let $\bar:\CC(v)@>>>\CC(v)$ be the field automorphism which is the identity on
$\CC$ and maps $v$ to $v\i$. For $(w,\l)\in\tX_m$ we write $E(w,\l)=(-1)^{|u|}$ where 

(a) $w=zu$, $(z,\l)\in\tX_m^0$, $u\in I_{z,\l}\sub W_\l$,
\nl
see 5.2(b).

We show:

(b) If $(w,\l)\in\tX_m,s\in S$ then $E(sws,s\l)=E(w,\l)$;

(c) if $(w,\l)\in\tX_m,s\in S$ are such that $sw=ws$ and $s\in W_\l$ then $E(ws,\l)=-E(w,\l)$.
\nl
We write $w=zu$ as in (a). Assume first that $s\in W_\l$. We have $sws=z\io_z(s)us$ and $\io_z(s)\in W_{z(\l)}=W_\l=W_{s\l}$
hence $\io_z(s)us\in I_{z,\l}$ and $E(sws,s\l)=(-1)^{|\io_z(s)us|}=(-1)^{|u|}=E(w,\l)$. If $sw=ws$, we have $ws=zus$ and $us\in I_{z,\l}$
hence $E(ws,\l)=(-1)^{|us|}=-(-1)^{|u|}=-E(w,\l)$. Next we assume that $s\n W_\l$ then $s\in[\l,\l]$ hence $(szs,s\l)\in\tX_m^0$.
Moreover, $sws=szus=(szs)(sus)$ and $sus\in W_{s\l}$ and more precisely $sus\in I_{szs,s\l}$. Hence
$E(sws,s\l)=(-1)^{|sus|}=(-1)^{|u|}=E(w,\l)$. This proves (b) and (c).

\mpb

Clearly, there is a unique $\CC$-linear map $B:\MM_m@>>>\MM_m$ such that for any $(w,\l)\in\tX_m$ and any $f\in\CC(v)$ we have
$$B(f\aa_{w,\l})=\bar f E(w,\l)\ct_w\i\aa_{w,-m\l}.$$
 We state the main result of this section.

\proclaim{Proposition 7.2} (a) For any $s\in S$ and any $\x\in\MM_m$ we have $B(\ct_s\x)=\ct_s\i B(x)$.

(b) The square of the map $\bar{}:\MM_m@>>>\MM_m$ is equal to $1$.
\endproclaim
To prove (a) it is enough to show that for any $(w,\l)\in\tX_m$ and any $s\in S$ we have

(c) $B(\ct_s\aa_{w,\l})=E(w,\l)\ct_s\i\ct_w\i\aa_{w,-m\l}$.
\nl
We set $\D=1$ if $s\in W_\l$ and $\D=0$ if $s\n W_\l$.

Assume that $sw\ne ws$, $|sw|>|w|$. We have
$$B(\ct_s\aa_{w,\l})=B(\aa_{sws,s\l})=E(sws,s\l)\ct_{sws}\i\aa_{sws,-ms\l},$$
$$\align&E(w,\l))\ct_s\i\ct_w\i\aa_{w,-m\l}=E(w,\l)\ct_s\i\ct_w\i\ct_s\i\ct_s\aa_{w,-m\l}\\&=E(sws,s\l)\ct_{sws}\i\aa_{sws,-ms\l}.\endalign$$
Hence (c) holds in this case.

Assume that $sw\ne ws$, $|sw|<|w|$. We must show that
$$B(\aa_{sws,s\l}+\D(v^2-v^{-2})\aa_{w,\l})=E(w,\l)\ct_s\i \ct_w\i\aa_{w,-m\l},$$
or that
$$\align&E(sws,s\l)\ct_{sws}\i\aa_{sws,-ms\l}+\D(v^{-2}-v^2)E(w,\l)\ct_w\i\aa_{w,-m\l}\\&=E(w,\l)\ct_s\i\ct_w\i\aa_{w,-m\l},\endalign$$
or that
$$\align&\ct_{sws}\i \aa_{sws,s\l}+\D(v^{-2}-v^2)\ct_s\i\ct_{sws}\i \ct_s\i\aa_{w,-m\l}\\&=  
\ct_s\i\ct_s\i\ct_{sws}\i\ct_s\i\aa_{w,-\l},\endalign$$
or that
$$\ct_s\ct_s\ct_{sws}\i\aa_{sws,-s\l}+\d(v^{-2}-v^2)\ct_s\ct_{sws}\i \ct_s\i\aa_{w,-m\l}=\ct_{sws}\i\ct_s\i\aa_{w,-m\l},$$
or that
$$\align&\ct_{sws}\i\aa_{sws,-ms\l}+(v^2-v^{-2})\D\ct_s\ct_{sws}\i\aa_{sws,-ms\l}
\\&+\D(v^{-2}-v^2)\ct_s\ct_{sws}\i \ct_s\i\aa_{w,-m\l}=\ct_{sws}\i\ct_s\i\aa_{w,-m\l}.\endalign$$
Here we substitute $\ct_{sws}\i\aa_{sws,-ms\l}=\ct_{sws}\i\ct_s\i\aa_{w,-m\l}$.
It remains to note that
$$\align&\ct_{sws}\i\ct_s\i\aa_{w,-m\l}+(v^2-v^{-2})\D\ct_s\ct_{sws}\i\ct_s\i\aa_{w,-m\l}\\&
+\D(v^{-2}-v^2)\ct_s\ct_{sws}\i\ct_s\i\aa_{w,-m\l}=\ct_{sws}\i\ct_s\i\aa_{w,-m\l}.\endalign$$
This proves (c) in our case.

Assume that $sw=ws$, $|sw|>|w|$. We must show that
$$B(\aa_{w,s\l}+\D(v+v\i)\aa_{sw,\l})=E(w,\l)\ct_s\i \ct_w\i\aa_{w,-m\l},$$
or that
$$\align&E(w,\l)\ct_w\i\aa_{w,-ms\l}+\D(v+v\i)E(sw,\l)\ct_{sw}\i\aa_{sw,-m\l}\\&=
E(w,\l)\ct_s\i \ct_w\i\aa_{w,-m\l},\endalign$$
or that
$$\align&\ct_w\i\aa_{w,-ms\l}-\D(v+v\i)\ct_w\i\ct_s\i\aa_{sw,-m\l}
\\&=\ct_s\i\ct_w\i\aa_{w,-m\l},\endalign$$
or that
$$\aa_{w,-ms\l}-\D(v+v\i)\ct_s\i\aa_{sw,-m\l}=\ct_s\i\aa_{w,-m\l},$$
or that
$$\ct_s\aa_{w,-ms\l}-\D(v+v\i)\aa_{sw,-m\l}=\aa_{w,-m\l},$$
or that
$$\ct_s\aa_{w,-ms\l}=\aa_{w,-m\l}+\D(v+v\i)\aa_{sw,-m\l}.$$
This follows from the definitions.
This proves (c) in our case.

Assume that $sw=ws$, $|sw|<|w|$. We must show that
$$\align&B(\D(v-v\i)\aa_{sw,\l}+\D(v^2-v^{-2}-1)\aa_{w,\l}+(1-\D)\aa_{w,s\l})\\&
=E(w,\l)\ct_s\i \ct_w\i\aa_{w,-m\l},\endalign$$
or that 
$$\align&\D(v\i-v)E(sw,\l)\ct_{sw}\i\aa_{sw,-m\l}+\D(v^{-2}-v^2-1)E(w,\l)\ct_w\i\aa_{w,-m\l}\\&+(1-\D)E(w,s\l)\ct_w\i\aa_{w,-ms\l})
=E(w,\l)\ct_s\i \ct_w\i\aa_{w,-m\l},\endalign$$
or that 
$$\align&\D(v\i-v)\ct_{sw}\i\aa_{sw,-m\l}-\D(v^{-2}-v^2-1)\ct_{sw}\i\ct_s\i\aa_{w,-m\l}\\&-(1-\D)\ct_{sw}\i\ct_s\i\aa_{w,-ms\l})
=-\ct_{sw}\i\ct_s\i\ct_s\i\aa_{w,-m\l},\endalign$$
or that 
$$\align&\D(v\i-v)\aa_{sw,\l}-\D(v^{-2}-v^2-1)\ct_s\i\aa_{w,-m\l}-(1-\D)\ct_s\i\aa_{w,-ms\l})\\&=-\ct_s\i\ct_s\i\aa_{w,-m\l}
\endalign$$
or that 
$$\D(v\i-v)\ct_s\ct_s\aa_{sw,-m\l}-\D(v^{-2}-v^2-1)\ct_s\aa_{w,-m\l}-(1-\D)\ct_s\aa_{w,-ms\l})=-\aa_{w,-m\l}$$
or that 
$$\align&\D(v\i-v)\aa_{sw,-m\l}+\D(v\i-v)(v^2-v^{-2}-1)\ct_s\aa_{sw,-m\l}\\&
-\D(v^{-2}-v^2-1)\ct_s\aa_{w,-m\l}-(1-\D)\ct_s\aa_{w,-ms\l})=-\aa_{w,-m\l}.\endalign$$
When $\D=0$ this is just $\ct_s\aa_{w,-ms\l})=\aa_{w,-m\l}$ which follows from the definitions. When $\D=1$ we see
that it is enough to observe the following obvious equality:
$$\align&(v\i-v)\aa_{sw,-m\l}+(v\i-v)(v^2-v^{-2})(\aa_{sw,-m\l}
+(v+v\i)\aa_{w,-m\l})\\&+(v^2-v^{-2}+1)((v-v\i)\aa_{sw,-m\l}
+(v^2-v^{-2}-1))\aa_{w,-m\l})=-\aa_{w,-m\l}.\endalign$$
This completes the proof of (c) hence that of (a).

We prove (b). We first show that for $(w,\l)\in\tX_m$ and $s\in S$ we have
$$B(\ct_s\i\aa_{w,\l})=\ct_s B(\aa_{w,\l}).\tag d$$
Indeed, the left hand side equals
$B(\ct_s\aa_{w,\l})+B((v^2-v^{-2})\aa_{w,\l})$ which by (a) equals 
$\ct_s\i B(\aa_{w,\l})+(v^{-2}-v^2)B(\aa_{w,\l})$ and this equals
$\ct_s B(\aa_{w,\l})$.
Using (d) repeatedly we see that $B(\ct_{w'}\i\aa_{w,\l})=\ct_{w'} B(\aa_{w,\l})$ for any $w'\in W$.
To prove (b) it is enough to prove that for any $(w,\l)\in\tX_m$ we have
$$B(B(\aa_{w,\l}))=\aa_{w,\l}$$
that is,
$$B(\ct_w\i\aa_{w,-m\l})=E(w,\l)\aa_{w,\l}.$$
The left hand side is equal to $\ct_w B(\aa_{w,-m\l})$ hence to
$$E(w,\l)\ct_w\ct_w\i\aa_{w,\l}=E(w,\l)\aa_{w,\l}.$$
 This completes the proof of (b).

\subhead 7.3\endsubhead
Let $(z,\l)\in\tX_m^0$. We show:

(a) $B(\aa_{z,\l})=\aa_{z,\l}$.
\nl
We must show that $\ct_z\i\aa_{z,-m\l}=\aa_{z,\l}$ 
or that $\ct_z\aa_{z,\l}=\aa_{z,-m\l}$. This follows the definition of the $\HH_m$-module structure
on $\MM_m$ since $zzz\i=z,z(\l)=-m\l$.

\subhead 7.4\endsubhead
Let $\cl$ be the $\ZZ[v\i]$-submodule of $\MM_m$ with basis 
$\{\aa_{w,\l};(w,\l)\in\tX_m\}$. From 7.2 one can deduce (a),(b) below by standard arguments.

(a) For any $(w,\l)\in\tX_m$ there is a unique element $\hat\aa_{w,\l}\in\MM_m$ such that

(i) $\hat\aa_{w,\l}\in\cl$, $\hat\aa_{w,\l}-\aa_{w,\l}\in v\i\ZZ[v\i]$,

(ii) $B(\hat\aa_{w,\l})=\hat\aa_{w,\l}$.
\nl
Moreover,

(b) $\{\hat\aa_{w,\l};(w,\l)\in\tX_m\}$ is a $\ZZ[v\i]$-basis of $\cl$ and a
$\CC(v)$-basis of $\MM_m$.
\nl
For example if $(z,\l)\in\tX_m^0$ then $\hat\aa_{z,\l}=\aa_{z,\l}$.

\enddocument